\documentclass[11pt,a4paper]{amsart}
\usepackage[left=25mm,top=30mm,right=25mm,bottom=20mm]{geometry}
\usepackage[utf8]{inputenc}

\usepackage{amssymb}
\usepackage[all]{xy}
\usepackage{indentfirst}
\usepackage{amsmath,amscd}
\usepackage{tikz,tikz-cd}

\allowdisplaybreaks[2]

\newtheorem{theorem}{Theorem}[section]
\newtheorem{proposition}[theorem]{Proposition}
\newtheorem{lemma}[theorem]{Lemma}
\newtheorem{corollary}[theorem]{Corollary}

\theoremstyle{definition}
\newtheorem{remark}[theorem]{Remark}
\newtheorem{definition}[theorem]{Definition}
\newtheorem{example}[theorem]{Example}
\numberwithin{equation}{section}

\begin{document}
\title{Cohomology of Leibniz Triple Systems and its applications}
\author{Xueru Wu}
\address{X. Wu: School of
Mathematics and Statistics, Northeast Normal
 University, Changchun 130024, China}
\email{wuxr884@nenu.edu.cn}
\author{Liangyun Chen}
\address{L. Chen: School of
Mathematics and Statistics, Northeast Normal
 University, Changchun 130024, China}
\email{chenly640@nenu.edu.cn}
\author{Yao Ma$^{*}$}
\address{Y. Ma: School of
Mathematics and Statistics, Northeast Normal
 University, Changchun 130024, China}
\email{may703@nenu.edu.cn}

\thanks{*Corresponding author.}

\thanks{\emph{MSC}(2020). 17A32, 17A40, 17B56, 16S20.}
\thanks{\emph{Key words and phrases}. Leibniz triple systems, cohomology, central extension, $T^*$-extension, deformation.}

\thanks{The second author is supported by NNSF of China (No.11771069). The third author is supported by NNSF of China (Nos. 11801066, 11771410).}

\begin{abstract}
In this paper, we introduce the first and third cohomology groups on Leibniz triple systems, which can be applied to extension theory and $1$-parameter formal deformation theory. Specifically, we investigate the central extension theory for Leibniz triple systems and show that there is a one-to-one correspondence between equivalent classes of central extensions of Leibniz triple systems and the third cohomology group. We study the $T^*$-extension of a Leibniz triple system and we determined that every even-dimensional quadratic Leibniz triple system $(\mathfrak{L},B)$ is isomorphic to a $T^*$-extension of a Leibniz triple system under a suitable condition. We also give a necessary and sufficient condition for a quadratic Leibniz triple system to admit a symplectic form. At last, we develop the $1$-parameter formal deformation theory of Leibniz triple systems and prove that it is governed by the cohomology groups.
\end{abstract}
\maketitle

\section{Introduction}

The notion of Leibniz triple systems was introduced by Bremner and S\'{a}nchez-Ortega {\rm\cite{Bremner}} in 2014, which were
obtained by using Kolesnikov-Pozhideav algorithm that takes the defining identities for a variety of algebras and produces the defining identities for the corresponding variety of dialgebras {\rm\cite{Co4}}, to Lie triple systems. Besides, Leibniz triple systems can also be obtained from any Malcev dialgebra {\rm\cite{BPS}}. Furthermore, Leibniz triple systems are related to Leibniz algebras in the same way that Lie triple
systems are related to Lie algebras. So it is natural to prove the analogous results from the theory of Lie triple systems to Leibniz triple systems. We list here some known results on this kind of algebraic structures: a Leibniz triple system decomposes as the orthogonal direct sum of well-described ideals if it admits a multiplicative basis \cite{Di};
Levi's theorem for Leibniz triple systems is determined \cite{Ma}; the theory of centroid on Leibniz triple systems refers to \cite{Cao}. In this paper, we will introduce the cohomology theory of Leibniz triple systems and use it to study the extension theory and $1$-parameter formal deformation theory, which generalize partial results in \cite{Ba,Bs,G1,LCC,Lin,Ma1,Sun,WZ,YB,Yam}.

%Cohomologies are useful tools in mathematic, which can be applied to algebras and to topologies as well as the theory of smooth manifolds or of holomorphic functions. The cohomology of Lie algebras was defined by Chevalley and Eilenberg in order to give an algebraic construction of the cohomology of the underlying topological spaces of compact Lie groups {\rm\cite{CE}}.
%The cohomology theory was then generalized to Lie superalgebras \cite{Scheunert}, Leibniz algebras \cite{Ad}, Lie triple systems \cite{Yam} and $n$-Lie algebras \cite{ANM}. This paper gives dijieshangtongdiao.

The paper is organized as follows. In Section 2, first, we will recall some basic definitions of Leibniz triple system and then we introduce the first and third cohomology groups of Leibniz triple systems. Section 3 concerns central extensions and we prove that there is a one-to-one correspondence between equivalent classes of central extensions of a Leibniz triple system and the third cohomology group. Section 4 is to introduce the $T^*$-extension theory of Leibniz triple systems. We prove that every even-dimensional quadratic Leibniz triple system containing an isotropic ideal is isomorphic to the $T^*$-extension of a Leibniz triple system. We also give the definition of symplectic form for a Leibniz triple system and get a sufficient and necessary condition for a quadratic Leibniz triple system to admit a symplectic structure. In Section 5, we study the $1$-parameter formal deformation theory of Leibniz triple systems. We show that the cohomology group defined in Section 2 is suitable for this $1$-parameter formal deformation theory.

Throughout this paper, all Leibniz triple system $\mathfrak{L}$ are defined over a fixed but arbitrary field $\mathbb{F}$.

\section{Cohomology of Leibniz triple systems}

In this section, we first recall some basic definitions and give the dual modules of Leibniz triple systems. Then we will give the first and third cohomology group of Leibniz triple systems.

\begin{definition}{\rm\cite{Bremner}}
A Leibniz triple system is a vector space $\mathfrak{L}$ endowed with a trilinear operation $ \{ \cdot, \cdot, \cdot \}: \mathfrak{L}\times \mathfrak{L}\times \mathfrak{L}\longrightarrow \mathfrak{L}$ satisfying
\begin{align}
& \{ a, b, \{ c, d, e \} \} = \{ \{ a, b, c \}, d, e \} - \{ \{ a, b, d \}, c, e \} - \{ \{ a, b, e \}, c, d \} + \{ \{ a, b, e \}, d, c \},\label{222.2}\\
& \{ a, \{ b, c, d \}, e \} = \{ \{ a, b, c \}, d, e \} - \{ \{ a, c, b \}, d, e \} - \{ \{ a, d, b \}, c, e \} + \{ \{ a, d, c \}, b, e \},\label{222.1}
\end{align}
for all $a, b, c, d, e \in \mathfrak{L}$.
\end{definition}

Note that a Leibniz triple system can be given by a Lie triple system with the same ternary product. A Leibniz algebra $L$ with product $[\cdot, \cdot]$ becomes a Leibniz triple system when $\{x, y, z\}:= [[x, y], z],$ for all $x, y, z\in L.$ More examples refer to {\rm\cite{Bremner}}. Denote by End$(\mathfrak{L})$ the set consisting of all linear maps on a Leibniz triple system $(\mathfrak{L}, \{\cdot, \cdot, \cdot\})$.

\begin{definition}{\rm\cite{Di}}
Let $I$ be a subspace of a Leibniz triple system $\mathfrak{L}.$ Then $I$ is called a subsystem of $\mathfrak{L},$ if $\{I,I,I\}\subset I;$ $I$ is called an ideal of $\mathfrak{L},$ if $\{I,\mathfrak{L},\mathfrak{L}\} + \{\mathfrak{L},I,\mathfrak{L}\} + \{\mathfrak{L},\mathfrak{L},I\}\subset I.$ Moreover, if $\{I, I, \mathfrak{L}\}=\{ I, \mathfrak{L},I\}=\{\mathfrak{L},I,I\}=0,$ then $I$ is called an abelian ideal of $\mathfrak{L}.$
\end{definition}

\begin{definition}{\rm\cite{Ma}}
Let $\mathfrak{L}$ be a Leibniz triple system and $V$ a vector space. $V$ is called an $\mathfrak{L}$-module, if $\mathfrak{L}\dot{+}V$ is a Leibniz triple system such that $(1)$ $\mathfrak{L}$ is a subsystem, $(2)$ $\{a,b,c\}\in V$ if any one of $a,b,c \in V;$ $(3)$ $\{a,b,c\}=0$ if any two of $a,b,c \in V.$
\end{definition}

\begin{definition}{\rm\cite{Ma}}
Let $\mathfrak{L}$ be a Leibniz triple system and $V$ a vector space. Suppose $l,m, r : \mathfrak{L} \times \mathfrak{L} \longrightarrow$ End$(V)$ are bilinear maps such that
\begin{align}
l(a,\{b, c, d\})&\!=\!l(\{a, b, c\},d)\! -\! l(\{a, c, b\},d)\! -\! l(\{a,d,b\},c)\! +\! l(\{a,d,c\},b),\label{222.11}\\
m(a, d)l(b, c)&=m(\{a, b, c\}, d) - m(\{a, c, b\}, d) - r(c, d)m(a, b) + r(b, d)m(a, c),\label{222.12}\\
m(a, d)m(b, c)&=r(c, d)l(a, b) - r(c, d)m(a, b) - m(\{a, c, b\}, d) + r(b, d)l(a, c),\label{222.13}\\
m(a, d)r(b, c)&=r(c, d)m(a, b) - r(c, d)l(a, b) - r(b, d)l(a, c) + m(\{a, c, b\}, d),\label{222.14}\\
r(\{a, b, c\}, d)&=r(c, d)r(a, b) - r(c, d)r(b, a) - r(b, d)r(c, a) + r(a, d)r(c, b),\label{222.15}\\
l(a, b)l(c, d)&=l(\{a, b, c\}, d) - l(\{a, b, d\}, c) - r(c, d)l(a, b) + r(d, c)l(a, b),\label{222.16}\\
l(a, b)m(c, d)&=m(\{a, b, c\}, d) - r(c, d)l(a, b) - l(\{a, b, d\}, c) + m(\{a, b, d\}, c),\label{222.17}\\
l(a, b)r(c, d)&=r(c, d)l(a, b) - m(\{a, b, c\}, d) - m(\{a, b, d\}, c) + l(\{a, b, d\}, c),\label{222.18}\\
m(a,\! \{b,\! c,\! d\})&\!=\!r(c,\! d)m(a,\! b)\! -\! r(b,\! d)m(a,\! c)\! -\! r(b,\! c)m(a,\! d)\! +\! r(c,\! b)m(a,\! d),\label{222.19}\\
r(a, \{b, c, d\})&=r(c, d)r(a, b) - r(b, d)r(a, c) - r(b, c)r(a, d) + r(c, b)r(a, d),\label{222.20}
\end{align}
for all $a, b, c, d \in \mathfrak{L}.$ Then $(r, m, l)$ is called a representation of $\mathfrak{L}$ on $V$.
\end{definition}

\begin{remark}\label{xz2.6}
In {\rm\cite{Ma}}, the authors proved that the notions of an $\mathfrak{L}$-module and a representation of $\mathfrak{L}$ are equivalent in the following sense. If $(r,m,l)$ is a representation of $\mathfrak{L}$ on $V,$ then $\mathfrak{L} \dot{+} V$ endowed with
\begin{align}
\{x + u, y + v, z + w\}_{\mathfrak{L} \dot{+} V} = \{x, y, z\} + (r(y, z)u + m(x,z)v + l(x, y)w)\label{222.123}
\end{align}
becomes a Leibniz triple system. If $V$ is an $\mathfrak{L}$-module then it induces a representation of $\mathfrak{L}$ on $V$ by setting $l(a,b)(v)=\{a,b,v\},$ $m(a,b)(v)=\{a,v,b\}$ and $r(a,b)(v)=\{v,a,b\},$ for any $a,b\in\mathfrak{L}$ and $v\in V.$ In particular, $\mathfrak{L}$ together with $l(a,b)(c)=\{a,b,c\},$ $m(a,b)(c)=\{a,c,b\}$ and $r(a,b)(c)=\{c,a,b\},$ for all $a,b,c\in\mathfrak{L},$ gives a representation of $\mathfrak{L}$ on $\mathfrak{L},$ called the adjoint $\mathfrak{L}$-module. Next, we give the example of the dual $\mathfrak{L}$-module.
\end{remark}

\begin{example}\label{ex222.6}
Let $(\mathfrak{L}, \{\cdot, \cdot, \cdot\})$ be a Leibniz triple system and $\mathfrak{L}^{*}$  the dual space of $\mathfrak{L}.$ Define the following bracket on $\mathfrak{L}\dot{+}\mathfrak{L}^{*}$ as
\begin{align}
\{a + f, b + g, c + h\}_{\mathfrak{L} \dot{+} \mathfrak{L}^{*}} = \{a, b, c\} + \{f, b, c\} + \{a, g, c\} + \{a, b, h\},
\end{align}
for all $a, b, c \in \mathfrak{L}$ and $f, g, h \in \mathfrak{L}^{*},$ where $\{f, b, c\}(n)=f(\{n, c, b\}),$ $\{a, g, c\}(n)=g(\{c, n, a\})$ and $\{a, b, h\}(n)=h(\{b, a, n\}),$ for all $n\in \mathfrak{L}.$ Then $(\mathfrak{L}\dot{+}\mathfrak{L}^{*},\{\cdot, \cdot, \cdot\}_{\mathfrak{L} \dot{+} \mathfrak{L}^{*}})$ is a Leibniz triple system, which can be checked straightforwardly, here we only check Eq. (\ref{222.2}) as an example. Note that, for all $a, b, c, d, e \in \mathfrak{L}$ and $f, g, h, s, t \in \mathfrak{L}^{*},$ we have
\begin{align*}
&\{a+ f, b+ g,\{ c+ h, d+s, e+t\}_{\mathfrak{L} \dot{+} \mathfrak{L}^{*}}\}_{\mathfrak{L} \dot{+} \mathfrak{L}^{*}}\notag\\
=&\{a+ f, b+ g,\{c,d,e\}+\{h,d,e\}+\{c,s,e\}+\{c,d,t\}\}_{\mathfrak{L} \dot{+} \mathfrak{L}^{*}}\notag\\
=&\{a,b,\{c,d,e\}\}+\{a,b,\{h,d,e\}+\{c,s,e\}+\{c,d,t\}\}+\{f,b,\{c,d,e\}\} + \{a,g,\{c,d,e\}\},
\end{align*}
\begin{align*}
&\{\{a+ f, b+ g, c+ h\}_{\mathfrak{L} \dot{+} \mathfrak{L}^{*}}, d+s, e+t\}_{\mathfrak{L} \dot{+} \mathfrak{L}^{*}}\notag\\
=& \{\{a, b, c\}+ \{a,b,h\}+\{f,b,c\}+\{a,g,c\}, d+s, e+t\}_{\mathfrak{L} \dot{+} \mathfrak{L}^{*}}\notag\\
=& \{\{a, b, c\}, d, e\}+ \{\{a, b, c\},d,t\}+\{\{a, b, c\},s,e\}+\{\{a,b,h\}+\{f,b,c\}+\{a,g,c\},d,e\},
\end{align*}
\begin{align*}
&\{\{a+ f, b+ g, d+ s\}_{\mathfrak{L} \dot{+} \mathfrak{L}^{*}}, c+h, e+t\}_{\mathfrak{L} \dot{+} \mathfrak{L}^{*}}\notag\\
=& \{\{a, b, d\}+ \{a,b,s\}+\{f,b,d\}+\{a,g,d\}, c+h, e+t\}_{\mathfrak{L} \dot{+} \mathfrak{L}^{*}}\notag\\
=& \{\{a, b, d\}, c, e\}+ \{\{a, b, d\},c,t\}+\{\{a, b, d\},h,e\}+\{\{a,b,s\}+\{f,b,d\}+\{a,g,d\},c,e\},
\end{align*}
\begin{align*}
&\{\{a+ f, b+ g, e+ t\}_{\mathfrak{L} \dot{+} \mathfrak{L}^{*}}, c+h, d+s\}_{\mathfrak{L} \dot{+} \mathfrak{L}^{*}}\notag\\
=& \{\{a, b, e\}+ \{a,b,t\}+\{f,b,e\}+\{a,g,e\}, c+h, d+s\}_{\mathfrak{L} \dot{+} \mathfrak{L}^{*}}\notag\\
=& \{\{a, b, e\}, c, d\}+ \{\{a, b, e\},c,s\}+\{\{a, b, e\},h,d\}+\{\{a,b,t\}+\{f,b,e\}+\{a,g,e\},c,d\},
\end{align*}
\begin{align*}
&\{\{a+ f, b+ g, e+ t\}_{\mathfrak{L} \dot{+} \mathfrak{L}^{*}}, d+s, c+h\}_{\mathfrak{L} \dot{+} \mathfrak{L}^{*}}\notag\\
=& \{\{a, b, e\}+ \{a,b,t\}+\{f,b,e\}+\{a,g,e\}, d+s, c+h\}_{\mathfrak{L} \dot{+} \mathfrak{L}^{*}}\notag\\
=& \{\{a, b, e\}, d, c\}+ \{\{a, b, e\},d,h\}+\{\{a, b, e\},s,c\}+\{\{a,b,t\}+\{f,b,e\}+\{a,g,e\},d,c\},
\end{align*}
If $\mathfrak{L} \dot{+} \mathfrak{L}^{*}$ becomes a Leibniz triple system, we must have
\begin{align*}
&\{a+ f, b+ g,\{ c+ h, d+s, e+t\}_{\mathfrak{L} \dot{+} \mathfrak{L}^{*}}\}_{\mathfrak{L} \dot{+} \mathfrak{L}^{*}}\\
=& \{\!\{a+ f,\! b+ g,\! c+ h\}_{\mathfrak{L} \dot{+} \mathfrak{L}^{*}\!}, \!d+s,\! e+t\}_{\mathfrak{L}\! \dot{+} \!\mathfrak{L}^{*}}
 - \{\!\{\!a+ f,\! b+ g,\! d+ s\}_{\mathfrak{L} \dot{+} \mathfrak{L}^{*}\!},\! c+h,\! e+t\}_{\mathfrak{L}\! \dot{+} \! \mathfrak{L}^{*}}\\
& - \{\!\{\!a+ f,\! b+ g,\! e+ t\}_{\mathfrak{L} \dot{+} \mathfrak{L}^{*}\!},\! c+h,\! d+s\}_{\mathfrak{L}\! \dot{+} \! \mathfrak{L}^{*}}
+ \{\!\{\!a+ f,\! b+ g,\! e+ t\}_{\mathfrak{L} \dot{+} \mathfrak{L}^{*}\!},\! d+s, \!c+h\}_{\mathfrak{L}\! \dot{+}\! \mathfrak{L}^{*}}.
\end{align*}
Since
\begin{align*}
\{f, b, \{c, d, e\}\}(n)=&f(\{n,\{c, d, e\},b\}),\\
\{\{f, b, c\},d, e\}(n)=&\{f,b, c\}(\{n,e,d\})=f(\{\{n,e,d\},c,b\}),\\
\{\{f,b, d\},c, e\}(n)=&\{f,b, d\}(\{n,e,c\})=f(\{\{n,e,c\},d,b\}),\\
\{\{f,b, e\},c, d\}(n)=&\{f,b, e\}(\{n,d,c\})=f(\{\{n,d,c\},e,b\}),\\
\{\{f,b, e\},d,c\}(n)=&\{f,b, e\}(\{n,c,d\})=f(\{\{n,c,d\},e,b\}),
\end{align*}
and
\begin{align*}
\{n,\{c, d, e\},b\}=\{\{n,c,d\},e,b\}-\{\{n,d,c\},e,b\}-\{\{n,e,c\},d,b\}+\{\{n,e,d\},c,b\},
\end{align*}
it follows that
\begin{align*}
\{f, b, \{c, d, e\}\}=\{\{f, b, c\},d, e\}-\{\{f,b, d\},c, e\}-\{\{f,b, e\},c, d\}+\{\{f,b, e\},d,c\}.
\end{align*}
Similarly, we have
\begin{align*}
\{a, g, \{c, d, e\}\}=&\{\{a, g, c\},d, e\}-\{\{a,g, d\},c, e\}-\{\{a,g, e\},c, d\}+\{\{a,g, e\},d,c\},\\
\{a, b, \{h, d, e\}\}=&\{\{a, b, h\},d, e\}-\{\{a,b, d\},h, e\}-\{\{a,b, e\},h, d\}+\{\{a,b, e\},d,h\},\\
\{a, b, \{c, s, e\}\}=&\{\{a, b, c\},s, e\}-\{\{a,b, s\},c, e\}-\{\{a,b, e\},s, c\}+\{\{a,b, e\},c,s\},\\
\{ a, b, \{c, d, t\}\}=&\{\{a, b, c\},d, t\}-\{\{a,b, d\},c, t\}-\{\{a,b, t\},c, d\}+\{\{a,b, t\},d,c\}.
\end{align*}
So we have proved that $\mathfrak{L} \dot{+} \mathfrak{L}^{*}$ satisfies Eq. (\ref{222.2}).

Similarly, $\mathfrak{L} \dot{+} \mathfrak{L}^{*}$ satisfies Eq. (\ref{222.1}), hence $(\mathfrak{L}\dot{+}\mathfrak{L}^{*},\{\cdot, \cdot, \cdot\}_{\mathfrak{L} \dot{+} \mathfrak{L}^{*}})$ is a Leibniz triple system, which implies that $\mathfrak{L}^{*}$ becomes an $\mathfrak{L}$-module, which is called the dual $\mathfrak{L}$-module.
\end{example}

\begin{definition}
Let $V$ be an $\mathfrak{L}$-module. A $(2n+1)$-linear map $f:\underbrace{\mathfrak{L}\times\cdots\times \mathfrak{L}}_{(2n+1)~{\rm times}}\rightarrow V$ is called a $(2n+1)$-cochain of $\mathfrak{L}$ on $V.$ Denote by $C^{2n+1}(\mathfrak{L},V)$ the set of all $(2n+1)$-cochains, for $n\geq 0.$
\end{definition}

\begin{definition}
Let $\mathfrak{L}$ be a Leibniz triple system and $(r,m,l)$ a representation of $\mathfrak{L}$ on $V.$ For $n=1,2$ the coboundary operator $\delta^{2n-1}:C^{2n-1}(\mathfrak{L},V)\rightarrow C^{2n+1}(\mathfrak{L},V)$ is defined as follows.\\
$(1)$ A $1$-coboundary operator of $\mathfrak{L}$ on $V$ is defined by
\begin{align*}
\delta^1:C^1(\mathfrak{L},V) &\rightarrow C^3(\mathfrak{L},V)\\
                            f&\mapsto \delta^1f
\end{align*}
for $f\in C^1(\mathfrak{L},V)$ and
\begin{align*}
\delta^1f(x_1,x_2,x_3)=r(x_2,x_3)f(x_1)+m(x_1,x_3)f(x_2)+l(x_1,x_2)f(x_3)-f(\{x_1,x_2,x_3\}).
\end{align*}
$(2)$ A $3$-coboundary operator of $\mathfrak{L}$ on $V$ consists a pair of maps $(\delta^3_1,\delta^3_2),$ where
\begin{align*}
\delta^3_i:C^3(\mathfrak{L},V)&\rightarrow C^5(\mathfrak{L},V)\\
                            g&\mapsto \delta^3_ig
\end{align*}
for $g\in C^3(\mathfrak{L},V)$ and
\begin{gather}
\begin{aligned}
&\delta^3_1g(x_1,x_2,x_3,x_4,x_5)\\
=&g(x_1,x_2,\{x_3,x_4,x_5\})-g(\{x_1,x_2,x_3\},x_4,x_5)+g(\{x_1,x_2,x_4\},x_3,x_5)\\
&+g(\{x_1,x_2,x_5\},x_3,x_4)-g(\{x_1,x_2,x_5\},x_4,x_3)+l(x_1,x_2)g(x_3,x_4,x_5)\\
&-r(x_4,x_5)g(x_1,x_2,x_3)+r(x_3,x_5)g(x_1,x_2,x_4)+r(x_3,x_4)g(x_1,x_2,x_5)\\
&-r(x_4,x_3)g(x_1,x_2,x_5),\label{222.221}
\end{aligned}
\end{gather}
\begin{gather}
\begin{aligned}
&\delta^3_2g(x_1,x_2,x_3,x_4,x_5)\\
=&g(x_1,\{x_2,x_3,x_4\},x_5)-g(\{x_1,x_2,x_3\},x_4,x_5)+g(\{x_1,x_3,x_2\},x_4,x_5)\\
&+g(\{x_1,x_4,x_2\},x_3,x_5)-g(\{x_1,x_4,x_3\},x_2,x_5)+m(x_1,x_5)g(x_2,x_3,x_4)\\
&-r(x_4,x_5)g(x_1,x_2,x_3)+r(x_4,x_5)g(x_1,x_3,x_2)+r(x_3,x_5)g(x_1,x_4,x_2)\\
&-r(x_2,x_5)g(x_1,x_4,x_3).\label{222.222}
\end{aligned}
\end{gather}
\end{definition}

\begin{remark}
Similarly to the definition of $(r,m,l)$ in Remark $\ref{xz2.6}$, we can rewrite Eqs. $(\ref{222.221})$-$(\ref{222.222})$ as follows
\begin{gather}
\begin{aligned}
&\delta^3_1g(x_1,x_2,x_3,x_4,x_5)\\
=&g(x_1,x_2,\{x_3,x_4,x_5\})-g(\{x_1,x_2,x_3\},x_4,x_5)+g(\{x_1,x_2,x_4\},x_3,x_5)\\
&+g(\{x_1,x_2,x_5\},x_3,x_4)-g(\{x_1,x_2,x_5\},x_4,x_3)+\{x_1,x_2,g(x_3,x_4,x_5)\}\\
&-\{g(x_1,x_2,x_3),x_4,x_5\}+\{g(x_1,x_2,x_4),x_3,x_5\}+\{g(x_1,x_2,x_5),x_3,x_4\}\\
&-\{g(x_1,x_2,x_5),x_4,x_3\},\label{222.21}
\end{aligned}
\end{gather}
\begin{gather}
\begin{aligned}
&\delta^3_2g(x_1,x_2,x_3,x_4,x_5)\\
=&g(x_1,\{x_2,x_3,x_4\},x_5)-g(\{x_1,x_2,x_3\},x_4,x_5)+g(\{x_1,x_3,x_2\},x_4,x_5)\\
&+g(\{x_1,x_4,x_2\},x_3,x_5)-g(\{x_1,x_4,x_3\},x_2,x_5)+\{x_1,g(x_2,x_3,x_4),x_5\}\\
&-\{g(x_1,x_2,x_3),x_4,x_5\}+\{g(x_1,x_3,x_2),x_4,x_5\}+\{g(x_1,x_4,x_2),x_3,x_5\}\\
&-\{g(x_1,x_4,x_3),x_2,x_5\}.\label{222.22}
\end{aligned}
\end{gather}
\end{remark}

Recall that, a derivation of a Leibniz triple system $\mathfrak{L}$ is a linear map $D$ $:$ $\mathfrak{L}$ $\longrightarrow$ $\mathfrak{L}$ satisfying
\begin{align}
&D ( \{ a, b, c \} ) = \{ D ( a ), b, c \} + \{ a, D ( b ), c \} + \{ a, b, D ( c ) \},
\end{align}
for all $a, b, c \in \mathfrak{L}$. Now we introduce the first and third cohomology groups of a Leibniz triple system. We will see that the derivation is a $1$-cocycle.

Let $\mathfrak{L}$ be a Leibniz triple system and $V$ an $\mathfrak{L}$-module. We call the set
\begin{align*}
Z^1(\mathfrak{L},V)=\{f\in C^1(\mathfrak{L},V)~|~\delta^1f=0\}
\end{align*}
the space of $1$-cocycles of $\mathfrak{L}$ on $V$.
We can see that for any $f\in Z^1(\mathfrak{L},\mathfrak{L}),$ the space of $1$-cocycles of $\mathfrak{L}$ on its adjoint module, $f$ is a derivation of $(\mathfrak{L},\{\cdot,\cdot,\cdot\}).$ The set
\begin{align*}
Z^3(\mathfrak{L},V)=\{g\in C^3(\mathfrak{L},V)~|~\delta^3_1g=\delta^3_2g=0\}
\end{align*}
is called the space of $3$-cocycles of $\mathfrak{L}$ on $V$. The set
\begin{align*}
B^3(\mathfrak{L},V)=\{\delta^1f~|~f\in C^1(\mathfrak{L},V)\}
\end{align*}
is called the space of $1$-coboundaries of $\mathfrak{L}$ on $V$.

\begin{example}
Let $\mathfrak{L}=\mathbb{K}\langle x,y\rangle$ be a Leibniz triple system with $\{x,y,y\}=\{y,y,y\}=x$ and $\mathfrak{L}^*$ its dual space with the dual basis $\{x^*,y^*\}.$ Define a trilinear map $\theta:\mathfrak{L}\times \mathfrak{L}\times \mathfrak{L}\rightarrow \mathfrak{L}^*$ by $\theta(x,y,y)=\theta(y,y,y)=x^*+y^*,$ the rest are $0.$ We can easily prove that $\theta$ is a $3$-cocycle of $\mathfrak{L}$.
\end{example}

\begin{proposition}
$Z^3(\mathfrak{L},V)$ and $B^3(\mathfrak{L},V)$ are defined as above. Then $B^3(\mathfrak{L},V)\subset Z^3(\mathfrak{L},V).$
\begin{proof}
For any $\delta^1f\in B^3(\mathfrak{L},V),$ we need to prove that $\delta^3_1\delta^1f=0$ and $\delta^3_2\delta^1f=0$. In fact, by the definitions of $1$-coboundary operator, $3$-coboundary operator, $\mathfrak{L}$-module and Eq. (\ref{222.21}), we have
\begin{align*}
&\delta^3_1\delta^1f(x_1,x_2,x_3,x_4,x_5)\\
=&\delta^1f(x_1,x_2,\{x_3,x_4,x_5\})-\delta^1f(\{x_1,x_2,x_3\},x_4,x_5)+\delta^1f(\{x_1,x_2,x_4\},x_3,x_5)\\
&+\delta^1f(\{x_1,x_2,x_5\},x_3,x_4)-\delta^1f(\{x_1,x_2,x_5\},x_4,x_3)+\{x_1,x_2,\delta^1f(x_3,x_4,x_5)\}\\
&-\{\delta^1f(x_1,x_2,x_3),x_4,x_5\}+\{\delta^1f(x_1,x_2,x_4),x_3,x_5\}+\{\delta^1f(x_1,x_2,x_5),x_3,x_4\}\\
&-\{\delta^1f(x_1,x_2,x_5),x_4,x_3\}\\
=&\{f(x_1),x_2,\{x_3,x_4,x_5\}\}+\{x_1,f(x_2),\{x_3,x_4,x_5\}\}+\{x_1,x_2,f(\{x_3,x_4,x_5\})\}\\
&-f(\{x_1,x_2,\{x_3,x_4,x_5\}\})-\{f(\{x_1,x_2,x_3\}),x_4,x_5\}-\{\{x_1,x_2,x_3\},f(x_4),x_5\}\\
&-\{\{x_1,x_2,x_3\},x_4,f(x_5)\}+f(\{\{x_1,x_2,x_3\},x_4,x_5\})+\{f(\{x_1,x_2,x_4\}),x_3,x_5\}\\
&+\{\{x_1,x_2,x_4\},f(x_3),x_5\}+\{\{x_1,x_2,x_4\},x_3,f(x_5)\} -f(\{\{x_1,x_2,x_4\},x_3,x_5\})\\
&+\{f(\{x_1,x_2,x_5\}),x_3,x_4\}+\{\{x_1,x_2,x_5\},f(x_3),x_4\}+\{\{x_1,x_2,x_5\},x_3,f(x_4)\}\\ &-f(\{\{x_1,x_2,x_5\},x_3,x_4\})-\{f(\{x_1,x_2,x_5\}),x_4,x_3\}-\{\{x_1,x_2,x_5\},f(x_4),x_3\}\\
&-\{\{x_1,x_2,x_5\},x_4,f(x_3)\}+f(\{\{x_1,x_2,x_5\},x_4,x_3\})+\{x_1,x_2,\{f(x_3),x_4,x_5\}\}\\
&+\{x_1,x_2,\{x_3,f(x_4),x_5\}\}+\{x_1,x_2,\{x_3,x_4,f(x_5)\}\}- \{x_1,x_2,f(\{x_3,x_4,x_5\})\}\\
&-\{\{f(x_1),x_2,x_3\},x_4,x_5\}-\{\{x_1,f(x_2),x_3\},x_4,x_5\}-\{\{x_1,x_2,f(x_3)\},x_4,x_5\}\\
&+\{f(\{x_1,x_2,x_3\}),x_4,x_5\}+\{\{f(x_1),x_2,x_4\},x_3,x_5\}+\{\{x_1,f(x_2),x_4\},x_3,x_5\}\\
&+\{\{x_1,x_2,f(x_4)\},x_3,x_5\}-\{f(\{x_1,x_2,x_4\}),x_3,x_5\}+\{\{f(x_1),x_2,x_5\},x_3,x_4\}\\
&+\{\{x_1,f(x_2),x_5\},x_3,x_4\}+\{\{x_1,x_2,f(x_5)\},x_3,x_4\}-\{f(\{x_1,x_2,x_5\}),x_3,x_4\}\\
&-\{\{f(x_1),x_2,x_5\},x_4,x_3\}-\{\{x_1,f(x_2),x_5\},x_4,x_3\}-\{\{x_1,x_2,f(x_5)\},x_4,x_3\}\\
&+\{f(\{x_1,x_2,x_5\}),x_4,x_3\}\\
=&0.
\end{align*}
On the other hand, it follows from the definitions of $1$-coboundary operator, $3$-coboundary operator, $\mathfrak{L}$-module and Eq. (\ref{222.22}) that
\begin{align*}
&\delta^3_2\delta^1f(x_1,x_2,x_3,x_4,x_5)\\
=&\delta^1f(x_1,\{x_2,x_3,x_4\},x_5)-\delta^1f(\{x_1,x_2,x_3\},x_4,x_5)+\delta^1f(\{x_1,x_3,x_2\},x_4,x_5)\\
&+\delta^1f(\{x_1,x_4,x_2\},x_3,x_5)-\delta^1f(\{x_1,x_4,x_3\},x_2,x_5)+\{x_1,\delta^1f(x_2,x_3,x_4),x_5\}\\
&-\{\delta^1f(x_1,x_2,x_3),x_4,x_5\}+\{\delta^1f(x_1,x_3,x_2),x_4,x_5\}+\{\delta^1f(x_1,x_4,x_2),x_3,x_5\}\\
&-\{\delta^1f(x_1,x_4,x_3),x_2,x_5\}\\
=&\{f(x_1),\{x_2,x_3,x_4\},x_5\}+\{x_1,f(\{x_2,x_3,x_4\}),x_5\}+\{x_1,\{x_2,x_3,x_4\},f(x_5)\}\\
&-f(\{x_1,\{x_2,x_3,x_4\},x_5\})-\{f(\{x_1,x_2,x_3\}),x_4,x_5\}-\{\{x_1,x_2,x_3\},f(x_4),x_5\}\\
&-\{\{x_1,x_2,x_3\},x_4,f(x_5)\}+f(\{\{x_1,x_2,x_3\},x_4,x_5\})+\{f(\{x_1,x_3,x_2\}),x_4,x_5\}\\
&+\{\{x_1,x_3,x_2\},f(x_4),x_5\}+\{\{x_1,x_3,x_2\},x_4,f(x_5)\} -f(\{\{x_1,x_3,x_2\},x_4,x_5\})\\
&+\{f(\{x_1,x_4,x_2\}),x_3,x_5\}+\{\{x_1,x_4,x_2\},f(x_3),x_5\}+\{\{x_1,x_4,x_2\},x_3,f(x_5)\}\\ &-f(\{\{x_1,x_4,x_2\},x_3,x_5\})-\{f(\{x_1,x_4,x_3\}),x_2,x_5\}-\{\{x_1,x_4,x_3\},f(x_2),x_5\}\\
&-\{\{x_1,x_4,x_3\},x_2,f(x_5)\}+f(\{\{x_1,x_4,x_3\},x_2,x_5\})+\{x_1,\{f(x_2),x_3,x_4\},x_5\}\\
&+\{x_1,\{x_2,f(x_3),x_4\},x_5\}+\{x_1,\{x_2,x_3,f(x_4)\},x_5\}- \{x_1,f(\{x_2,x_3,x_4\}),x_5\}\\
&-\{\{f(x_1),x_2,x_3\},x_4,x_5\}-\{\{x_1,f(x_2),x_3\},x_4,x_5\}-\{\{x_1,x_2,f(x_3)\},x_4,x_5\}\\
&+\{f(\{x_1,x_2,x_3\}),x_4,x_5\}+\{\{f(x_1),x_3,x_2\},x_4,x_5\}+\{\{x_1,f(x_3),x_2\},x_4,x_5\}\\
&+\{\{x_1,x_3,f(x_2)\},x_4,x_5\}-\{f(\{x_1,x_3,x_2\}),x_4,x_5\}+\{\{f(x_1),x_4,x_2\},x_3,x_5\}\\
&+\{\{x_1,f(x_4),x_2\},x_3,x_5\}+\{\{x_1,x_4,f(x_2)\},x_3,x_5\}-\{f(\{x_1,x_4,x_2\}),x_3,x_5\}\\
&-\{\{f(x_1),x_4,x_3\},x_2,x_5\}-\{\{x_1,f(x_4),x_3\},x_2,x_5\}-\{\{x_1,x_4,f(x_3)\},x_2,x_5\}\\
&+\{f(\{x_1,x_4,x_3\}),x_2,x_5\}\\
=&0.
\end{align*}
Thus the proposition has been proved.
\end{proof}
\end{proposition}

Hence we can define $1$-cohomology space and $3$-cohomology space of $\mathfrak{L}$ as the factor spaces
\begin{align*}
H^1(\mathfrak{L},V):=&Z^1(\mathfrak{L},V).\\
H^3(\mathfrak{L},V):=&Z^3(\mathfrak{L},V)/B^3(\mathfrak{L},V).
\end{align*}

\begin{remark}\label{rm2.13}
Recall that $\mathfrak{L}\dot{+}V$ endowed with Eq. $(\ref{222.123})$ becomes a Leibniz triple system. Using the $Z^3(\mathfrak{L},V)$ defined above, we have a new Leibniz triple system structure on $\mathfrak{L}\dot{+}V$ by a $3$-cocycle.

For any $\theta\in Z^3(\mathfrak{L},V),$ set
\begin{align*}
\{(x,a),(y,b),(z,c)\}_\theta=(\{x,y,z\},\theta(x,y,z)).
\end{align*}
Note that
\begin{align*}
&\{(x,a),(y,b),\{(u,c),(v,d),(w,e)\}_\theta\}_\theta\\
=&\{(x,a),(y,b),(\{u,v,w\},\theta(u,v,w))\}_\theta=(\{x,y,\{u,v,w\}\},\theta(x,y,\{u,v,w\})\\
=&\Big(\{\{x,y,u\},v,w\},\theta(\{x,y,u\},v,w)\Big)-\Big(\{\{x,y,v\},u,w\},\theta(\{x,y,v\},u,w)\Big)\\
&-\Big(\{\{x,y,w\},u,v\},\theta(\{x,y,w\},u,v)\Big)+\Big(\{\{x,y,w\},v,u\},\theta(\{x,y,w\},v,u)\Big)\\
=&\{\{(x,a),(y,b),(u,c)\}_\theta,(v,d),(w,e)\}_\theta-\{\{(x,a),(y,b),(v,d)\}_\theta,(u,c),(w,e)\}_\theta\\
&-\{\{(x,a),(y,b),(w,e)\}_\theta,(u,c),(v,d)\}_\theta+\{\{(x,a),(y,b),(w,e)\}_\theta,(v,d),(u,c)\}_\theta,
\end{align*}
and
\begin{align*}
&\{(x,a),\{(y,b),(u,c),(v,d)\}_\theta,(w,e)\}_\theta\\
=&\{(x,a),(\{y,u,v\},\theta(y,u,v)),(w,e)\}_\theta=(\{x,\{y,u,v\},w\},\theta(x,\{y,u,v\},w)\\
=&\Big(\{\{x,y,u\},v,w\},\theta(\{x,y,u\},v,w)\Big)-\Big(\{\{x,u,y\},v,w\},\theta(\{x,u,y\},v,w)\Big)\\
&-\Big(\{\{x,v,y\},u,w\},\theta(\{x,v,y\},u,w)\Big)+\Big(\{\{x,v,u\},y,w\},\theta(\{x,v,u\},y,w)\Big)\\
=&\{\{(x,a),(y,b),(u,c)\}_\theta,(v,d),(w,e)\}_\theta-\{\{(x,a),(u,c),(y,b)\}_\theta,(v,d),(w,e)\}_\theta\\
&-\{\{(x,a),(v,d),(y,b)\}_\theta,(u,c),(w,e)\}_\theta+\{\{(x,a),(v,d),(u,c)\}_\theta,(y,b),(w,e)\}_\theta.
\end{align*}
Then, $(\mathfrak{L}_\theta,\{\cdot,\cdot,\cdot\}_\theta)$ is a Leibniz triple system.
\end{remark}

In the following sections, we will apply the cohomology theory of Leibniz triple systems to discuss the central extension, $T^*$-extension and $1$-parameter formal deformation of a Leibniz triple system.

\section{Central extension of Leibniz triple system}
In this section, we give the definition of central extension of a Leibniz triple system  $(\mathfrak{L},\{\cdot,\cdot,\cdot\})$ by an $\mathfrak{L}$-module $V$ and prove that there is a one-to-one correspondence between equivalent classes of central extensions of $(\mathfrak{L},\{\cdot,\cdot,\cdot\})$ and $H^3(\mathfrak{L},V).$

\begin{definition}
Suppose that $(\mathfrak{L},\{\cdot,\cdot,\cdot\})$ is a Leibniz triple system and $V$ is an $\mathfrak{L}$-module. $(\mathfrak{L}_c,\{\cdot,\cdot,\cdot\}_c)$ is called a central extension of $(\mathfrak{L},\{\cdot,\cdot,\cdot\})$ by $V$ if there is an exact row of Leibniz triple systems
$$\xymatrix
{0\ar[r]& V\ar[r]^\iota& \mathfrak{L}_c\ar[r]^\pi& \mathfrak{L}\ar[r]& 0}$$
such that $\iota(V)$ is contained in the center of $\mathfrak{L}_c,$ which is defined as
\begin{align*}
Z(\mathfrak{L}_c)=\{x\in
\mathfrak{L}_c~|~\{x,\mathfrak{L}_c,\mathfrak{L}_c\}_c=\{\mathfrak{L}_c,x,\mathfrak{L}_c\}_c=\{\mathfrak{L}_c,\mathfrak{L}_c,x\}_c=0\}.
\end{align*}
\end{definition}

\begin{definition}
If $(\mathfrak{L}_c^\prime,\{\cdot,\cdot,\cdot\}_c^\prime)$ is another central extension of $(\mathfrak{L},\{\cdot,\cdot,\cdot\})$ by $V$ such that there is an isomorphism $\phi:(\mathfrak{L}_c,\{\cdot,\cdot,\cdot\}_c)\rightarrow (\mathfrak{L}_c^\prime,\{\cdot,\cdot,\cdot\}_c^\prime),$ and the diagram
\begin{equation*}
\begin{CD}
0@>>>V@>\iota>> \mathfrak{L}_c @>\pi>>\mathfrak{L}@>>>0\\
@.@|@VV\phi V@|@.\\
0@>>>V\ @>\iota^\prime>>\mathfrak{L}_c^\prime@>\pi^\prime>> \mathfrak{L}@>>> 0
\end{CD}
\end{equation*}
commutes, then the two central extensions $(\mathfrak{L}_c,\{\cdot,\cdot,\cdot\}_c)$ and $(\mathfrak{L}_c^\prime,\{\cdot,\cdot,\cdot\}_c^\prime)$ are said to be equivalent.
\end{definition}

\begin{theorem}
There is a one-to-one correspondence between equivalent classes of central extensions of $(\mathfrak{L},\{\cdot,\cdot,\cdot\})$ by $V$ and $H^3(\mathfrak{L},V).$
\begin{proof}

{\it Step 1, construct a map  $\sigma_1$ from  Ext$_c(\mathfrak{L},V)$ the set of equivalent classes of central extensions of $\mathfrak{L}~by~V$ to $H^3(\mathfrak{L},V).$}

We first show that there exists a map from the set of central extensions of $(\mathfrak{L},\{\cdot,\cdot,\cdot\})$ by $V$ to $Z^3(\mathfrak{L},V).$ Suppose that $(\mathfrak{L}_c,\{\cdot,\cdot,\cdot\}_c)$ is a central extension of $(\mathfrak{L},\{\cdot,\cdot,\cdot\})$ by $V.$ Then we have the following exact sequence
$$\xymatrix
{0\ar[r]& V\ar[r]^\iota& \mathfrak{L}_c\ar[r]^\pi& \mathfrak{L}\ar[r]& 0}.$$

Since $\pi$ is surjective, there exists a linear map $s:\mathfrak{L}\rightarrow \mathfrak{L}_c,$ such that $\pi s=\rm{id}_\mathfrak{L}.$ Then for any $x,y,z\in \mathfrak{L},$
\begin{align*}
\pi(\{s(x),s(y),s(z)\}_c)-\pi s(\{x,y,z\})=\{\pi s(x),\pi s(y),\pi s(z)\}-\{x,y,z\}=0.
\end{align*}
It follows that $\{s(x),s(y),s(z)\}_c-s(\{x,y,z\})\in {\rm Ker} \pi=\iota(V)\subset Z(\mathfrak{L}_c),$ hence, for all $x_1,x_2,x_3,x_4,$ $x_5\in \mathfrak{L},$
\begin{align*}
\{\{s(x_1),s(x_2),s(x_3)\}_c,x_4,x_5\}_c&=\{s(\{x_1,x_2,x_3\}),x_4,x_5\}_c,\\
\{x_4,\{s(x_1),s(x_2),s(x_3)\}_c,x_5\}_c&=\{x_4,s(\{x_1,x_2,x_3\}),x_5\}_c,\\
\{x_4,x_5,\{s(x_1),s(x_2),s(x_3)\}_c\}_c&=\{x_4,x_5,s(\{x_1,x_2,x_3\})\}_c.
\end{align*}
Consider a map $g:\mathfrak{L}\times \mathfrak{L}\times \mathfrak{L}\rightarrow V,$
\begin{align*}
\iota g(x,y,z)=\{s(x),s(y),s(z)\}_c-s(\{x,y,z\}),
\end{align*}
which is well-defined, since $\iota$ is injective. Clearly, $g\in C^3(\mathfrak{L},V).$ Moreover, it is straightforwardly to check $\delta^3_1g=\delta^3_2g=0,$ and so $g\in Z^3(\mathfrak{L},V).$ Here we only check that $\delta^3_1g=0$ as an example. Note that it is sufficient to show that $\iota\delta^3_1g=0,$ since $\iota$ is injective. In fact,
\begin{align*}
&\iota\delta^3_1g(x_1,x_2,x_3,x_4,x_5)\\
=&\iota\Big(g(\{x_1,x_2,x_3\},x_4,x_5)-g(\{x_1,x_2,x_4\},x_3,x_5)-g(\{x_1,x_2,x_5\},x_3,x_4)\\
&+g(\{x_1,x_2,x_5\},x_4,x_3)-g(x_1,x_2,\{x_3,x_4,x_5\})\Big)\\
=&\{s(\{x_1,x_2,x_3\}),s(x_4),s(x_5)\}_c-s(\{\{x_1,x_2,x_3\},x_4,x_5\})\\
&-\{s(\{x_1,x_2,x_4\}),s(x_3),s(x_5)\}_c+s(\{\{x_1,x_2,x_4\},x_3,x_5\})\\
&-\{s(\{x_1,x_2,x_5\}),s(x_3),s(x_4)\}_c+s(\{\{x_1,x_2,x_5\},x_3,x_4\})\\
&+\{s(\{x_1,x_2,x_5\}),s(x_4),s(x_3)\}_c-s(\{\{x_1,x_2,x_5\},x_4,x_3\})\\
&-\{s(x_1),s(x_2),s(\{x_3,x_4,x_5\}\}_c+s(\{x_1,x_2,\{x_3,x_4,x_5\}\})\\
=&\{\{s(x_1),s(x_2),s(x_3)\}_c,s(x_4),s(x_5)\}_c-\{\{s(x_1),s(x_2),s(x_4)\}_c,s(x_3),s(x_5)\}_c\\
&-\{\{s(x_1),s(x_2),s(x_5)\}_c,s(x_3),s(x_4)\}_c+\{\{s(x_1),s(x_2),s(x_5)\}_c,s(x_4),s(x_3)\}_c\\
&-\{s(x_1),s(x_2),\{s(x_3),s(x_4),s(x_5)\}_c\}_c=0.
\end{align*}

We claim that the above arguments induce a map from Ext$_c(\mathfrak{L},V)$ to $H^3(\mathfrak{L},V).$ Indeed, suppose that $(\mathfrak{L}_c,\{\cdot,\cdot,\cdot\}_c)$ and $(\mathfrak{L}_c^\prime,\{\cdot,\cdot,\cdot\}_c^\prime)$ are equivalent central extensions of $(\mathfrak{L},\{\cdot,\cdot,\cdot\})$ by $V.$ Then we have the following diagram
$$\xymatrix
{0\ar[r]& V\ar[r]^\iota\ar@{=}[d]& \mathfrak{L}_c \ar[r]^\pi\ar[d]^{\phi}& \mathfrak{L} \ar[r]\ar@{=}[d]\ar@/^/[l]^{s}& 0\\
0\ar[r]& V\ar[r]^{\iota^\prime}& \mathfrak{L}_c^\prime \ar[r]^{\pi^\prime}& \mathfrak{L}\ar[r]\ar@/^/[l]^{s^\prime}& 0}$$
such that $\phi \iota=\iota^\prime$ and $\pi=\pi^\prime \phi,$ with $\phi$ being an isomorphism and $\pi s=\pi^\prime s^\prime=\rm{id}_\mathfrak{L}.$ Let $g,g^\prime$ be their corresponding $3$-cocycles constructed as above, i.e.,
\begin{align*}
\iota g(x,y,z)=&\{s(x),s(y),s(z)\}_c-s(\{x,y,z\}),\\
\iota^\prime g^\prime(x,y,z)=&\{s^\prime(x),s^\prime(y),s^\prime(z)\}_c-s^\prime(\{x,y,z\}).
\end{align*}
Then
\begin{align}
\iota^\prime g(x,y,z)=&\phi \iota g(x,y,z)=\phi(\{s(x),s(y),s(z)\}_c)-\phi s(\{x,y,z\}).\label{**}
\end{align}

Now define $f\in C^1(\mathfrak{L},V)$ by $$\iota^\prime f(x)=s^\prime(x)-\phi s(x),~\forall x\in \mathfrak{L},$$ which is reasonable since Im$\iota^\prime=$Ker$\pi^\prime$ and
\begin{align*}
\pi^\prime s^\prime (x)-\pi^\prime \phi s(x)=x-\pi s(x)=0,
\end{align*}
Since $s^\prime(x)-\phi s(x)=\iota^\prime f(x)\in Z(\mathfrak{L}_c^\prime),$
\begin{align}
\{s^\prime(x),s^\prime(y),s^\prime(z)\}_c^\prime=\{\phi s(x),\phi s(y),\phi s(z)\}_c^\prime= \phi(\{s(x),s(y),s(z)\}_c).\label{****}
\end{align}
Combing Eqs. (\ref{**}) and (\ref{****}),
\begin{align*}
\iota^\prime(g^\prime-g)(x,y,z)=-\iota^\prime f(\{x,y,z\})=\iota^\prime(\delta^1f)(x,y,z).
\end{align*}
So $g^\prime-g=\delta^1f\in B^3(\mathfrak{L},V).$ Hence we give a map from Ext$_c(\mathfrak{L},V)$ to $H^3(\mathfrak{L},V)$ denote by $\sigma_1.$

\emph{Step 2, construct a map  $\sigma_2$ from $H^3(\mathfrak{L},V)$ to $Ext_c(\mathfrak{L},V).$}

Let $g\in Z^3(\mathfrak{L},V)$ and $\mathfrak{L}_c=\mathfrak{L}+ V$ with
\begin{align*}
\{(x,a),(y,b),(z,c)\}_c=(\{x,y,z\},g(x,y,z)).
\end{align*}
 Then Remark \ref{rm2.13} says that $(\mathfrak{L}_c,\{\cdot,\cdot,\cdot\}_c)$ is a Leibniz triple system. Define $\iota:V\rightarrow \mathfrak{L}_c$ by $\iota(a)=(0,a)$ and $\pi:\mathfrak{L}_c\rightarrow \mathfrak{L}$ by $\pi(x,a)=x.$ Clearly, $(\mathfrak{L}_c,\{\cdot,\cdot,\cdot\}_c)$ is a central extension of $(\mathfrak{L},\{\cdot,\cdot,\cdot\})$ by $V.$ Then we have construct a map from $Z^3(\mathfrak{L},V)$ to Ext$_c(\mathfrak{L},V).$

Suppose $g,g^\prime \in Z^3(\mathfrak{L},V)$ and $g-g^\prime\in B^3(\mathfrak{L},V)$, i.e., there exist $f\in C^1(\mathfrak{L},V)$ satisfying $(g^\prime-g)=\delta^1f.$ Then $(g^\prime-g)(x,y,z)=-f(\{x,y,z\}).$ Let $(\mathfrak{L}_c,\{\cdot,,\cdot,\cdot\}_c)$ and$(\mathfrak{L}_c^\prime,\{\cdot,,\cdot,\cdot\}_c^\prime)$ be two central extensions of $(\mathfrak{L},\{\cdot,\cdot,\cdot\})$ by $V$ defined as above with respect to $g$ and $g^\prime,$ respectively. Then $\iota(a)=\iota^\prime(a)=(0,a)$ and $\pi(x,a)=\pi^\prime(x,a)=x.$ Consider a linear map
\begin{align*}
\phi:(\mathfrak{L}_c,\{\cdot,,\cdot,\cdot\}_c)&\rightarrow (\mathfrak{L}_c^\prime,\{\cdot,,\cdot,\cdot\}_c^\prime)\\
                                        (x,a)&\mapsto (x,a-f(x)).
\end{align*}
Then $\phi \iota(a)=\iota^\prime(a)$ and $\pi^\prime \phi (x,a)=\pi^\prime(x,a-f(x))=x=\pi(x,a).$ Next, we prove that $\phi$ is an isomorphism.

If $\phi (x,a)=\phi(\tilde{x},\tilde{a}),$ it follows that $(x,a-f(x))=(\tilde{x},\tilde{a}-f(\tilde{x})),$ that is, $x=\tilde{x}$ and $a-f(x)=\tilde{a}-f(\tilde{x}),$ then $a=\tilde{a},$ and so $\phi$ is injective; $\phi$ is obviously surjective. Note that
\begin{align*}
&\phi(\{(x,a),(y,b),(z,c)\}_c)\\
=&\phi(\{x,y,z\},g(x,y,z))=(\{x,y,z\},g(x,y,z)-f(\{x,y,z\}))\\
                              =&(\{x,y,z\},g^\prime(x,y,z)) =\{(x,a-f(x)),(y,b-f(y)),(z,c-f(z))\}_c^\prime\\
                             =&\{\phi(x,a),\phi(y,b),\phi(z,c)\}_c^\prime.
\end{align*}
Therefore, $(\mathfrak{L}_c,\{\cdot,\cdot,\cdot\}_c)$ and $(\mathfrak{L}_c^\prime,\{\cdot,\cdot,\cdot\}_c^\prime)$ are equivalent central extensions of $(\mathfrak{L},\{\cdot,\cdot,\cdot\})$ by $V.$ Hence we find a map from $H^3(\mathfrak{L},V)$ to Ext$_c(\mathfrak{L},V)$ by $\sigma_2.$

It is obvious that $\sigma_1\sigma_2={\rm id}_{H^3(\mathfrak{L},V)}$ and $\sigma_2\sigma_1={\rm id}_{{\rm Ext}_c(\mathfrak{L},V)}.$ We get the assertion of the theorem.
\end{proof}
\end{theorem}

\section{$T^*$-extension of Leibniz triple systems}
At the beginning of this section we will study the quadratic structure of Leibniz triple systems, and then we introduce the theory of $T^*$-extension to Leibniz triple systems, at last, we study the symplectic form of a Leibniz triple system.

\subsection{Quadratic Leibniz triple systems}
In this subsection, we will introduce the quadratic Leibniz triple system $(\mathfrak{L},B)$ and prove that $\mathfrak{L}$ can be decomposed into the direct sum of non-degenerate ideals. We also give the relation between central descending series $(C^{n}(\mathfrak{L}))_{n\geq 0}$ and central ascending series $(C_{n}(\mathfrak{L}))_{n\geq 0}.$

\begin{definition}
Let $\mathfrak{L}$ be a Leibniz triple system. A bilinear form $B$ on $\mathfrak{L}$ is said to be symmetric if $B(x, y) = B(y, x)$, for all $x, y \in \mathfrak{L}$. $B$ is said to be non-degenerate if $\mathfrak{L}^{\bot} = 0,$ where $\mathfrak{L}^{\bot} = \{x \in \mathfrak{L} ~|~ B(x, y) = B(y, x) =0, \forall y \in \mathfrak{L}\}.$
\end{definition}

For any $x,y\in \mathfrak{L},$ denote by $R_{(x,y)},$ $M_{(x,y)}$ and $L_{(x,y)}$ the right, middle and left multiplication operator on $\mathfrak{L},$ respectively, i.e., $R_{(x,y)},~M_{(x,y)},~L_{(x,y)}: \mathfrak{L}\times \mathfrak{L}\rightarrow$ End$(\mathfrak{L})$ are defined as $R_{(x, y)}z:=\{z, x, y\},$ $M_{(x, y)}z:=\{x, z, y\}$ and $L_{(x, y)}z:=\{x, y, z\},$ which can be used to characterize the invariance of $B.$

\begin{definition}
A symmetric bilinear form $B$ on a Leibniz triple system $\mathfrak{L}$ is said to be right-invariant $($resp. left-invariant or middle-invariant $)$ if
\begin{align*}
B(R_{(a,b)}x, y)&=B(x, R_{(b,a)}y)\\
 \Big({\rm resp}. ~B(L_{(a,b)}x, y)=B(x, L_{(b,a)}y), ~&{\rm or}
 ~B(M_{(a,b)}x, y)=B(x, M_{(b,a)}y) \Big),
\end{align*}
for all $x, y, a, b \in \mathfrak{L}$.
\end{definition}

\begin{remark}
If $B$ is symmetric then any two invariance of $B$ lead to the third invariance. Here, we only check that the middle-invariant can be deduced by the right-invariant and left-invariant as an example, which follows from for all $a, b, x, y \in \mathfrak{L},$
\begin{align*}
B(M_{(a,b)}x, y) =& B(L_{(a,x)}b, y)= B(b, L_{(x,a)}y)=B(b, R_{(a, y)}x)\\
=& B(R_{(y, a)}b, x)=B(M_{(b, a)}y,x)= B(x, M_{(b, a)}y).
\end{align*}
\end{remark}

\begin{definition}
Let $\mathfrak{L}$ be a Leibniz triple system together with a bilinear form $B$. We say that $(\mathfrak{L}, B)$ is a quadratic Leibniz triple system if $B$ is non-degenerate symmetric invariant. An ideal $I$ of $\mathfrak{L}$ is called non-degenerate if $B|_{I\times I}$ is non-degenerate. A quadratic Leibniz triple system $(\mathfrak{L}, B)$ is said to be $B$-irreducible if $\mathfrak{L}$ contains no nontrivial non-degenerate ideal.
\end{definition}

\begin{lemma}\label{213.6}
Let $( \mathfrak{L}, B)$ be a quadratic Leibniz triple system.

$(1)$ If $\mathfrak{L}$ is finite dimensional and $I$ is a subspace of $\mathfrak{L},$ then $\mathfrak{L} = I \dot{+} I^{\bot},$ $(I^{\bot})^{\bot}=I.$

$(2)$ If $I$ is an ideal of $\mathfrak{L},$ then $I^{\bot}$ is an ideal of $\mathfrak{L}.$

$(3)$ If $\mathfrak{L}$ is finite dimensional and $I$ is a non-degenerate subspace of $\mathfrak{L},$ then $I^{\bot}$ is non-degenerate.
\begin{proof}
(1)
Let $\{x_1,x_2, ..., x_l\}$ be a basis of $I$ and extend it to a basis of $\mathfrak{L}$ say $\{x_1,x_2, ..., x_n\}.$ Then $x\in I^{\bot},$ if and only if  $B(x_i,x)=0,$ $1\leq i\leq l.$ Suppose $x=\sum\limits_{i=1}^{n}k_ix_i.$ Then $B(x_i,x)=0,$ if and only if $B(x_i,\sum\limits_{j=1}^{n}k_jx_j)=\sum\limits_{j=1}^{n}B(x_i,x_j)k_j=0,$ $1\leq i\leq l,$ i.e.,
\begin{equation*}       %开始数学环境
\left(                 %左括号
  \begin{array}{cccc}   %该矩阵一共3列，每一列都居中放置
    B(x_1,x_1) & \cdots &B(x_1,x_n)\\  %第一行元素
    \vdots&&\vdots\\
   B(x_l,x_1) &\cdots &B(x_l,x_n)\\  %第二行元素
  \end{array}
\right)                 %右括号
\left(                 %左括号
  \begin{array}{cccc}   %该矩阵一共3列，每一列都居中放置
    k_1 \\  %第一行元素
    \vdots\\
   k_n\\  %第二行元素
  \end{array}
\right)=\left(                 %左括号
  \begin{array}{cccc}   %该矩阵一共3列，每一列都居中放置
    0 \\  %第一行元素
    \vdots\\
   0\\  %第二行元素
  \end{array}
\right),                 %右括号
\end{equation*}
which implies that  $I^\bot$ is isomorphism to $\{X\in\mathbb{F}^n~|~AX=0\}$ as vector space. Then dim$I^\bot=n-$rank$A.$ Note that rank$A=l=$dim$I$  and $I\cap I^\bot=0$ by non-degenerate of $B.$ It follows that $\mathfrak{L} = I \dot{+} I^{\bot}.$
$(I^{\bot})^{\bot}=I$ follows from the facts that $I\subset(I^\bot)^\bot$ and dim$I=$dim$(I^\bot)^\bot=$dim$\mathfrak{L}-$dim$I^\bot.$

(2) It follows the invariance of $B.$

(3) Follows from the non-degenerate of $B$ and $B|_{I\times I}$ and $B(I,I^\bot)=0.$
\end{proof}
\end{lemma}

\begin{proposition}
Let $(\mathfrak{L}, B)$ be a finite dimensional quadratic Leibniz triple system. Then $\mathfrak{L} = \oplus_{i=1}^{r}\mathfrak{L}_{i}$ such that for all $1\leq i\leq r,$ we have

$(1)$ $\mathfrak{L}_{i}$ is a non-degenerate ideal$;$

$(2)$ $\mathfrak{L}_{i}$ is $B$-irreducible.
\begin{proof}
By induction on dim$\mathfrak{L}=n$ and using the facts that $\mathfrak{L} = I \dot{+} I^{\bot}$ and any ideal of $I$ is also an ideal of $\mathfrak{L},$ the proposition follows.

If $n=1,$ the assertion is true.

Assume that dim$\mathfrak{L}<n$ is true.

When dim$\mathfrak{L}= n,$ we will do it in two cases, if $\mathfrak{L}$ is $B$-irreducible, then the assertion is true; if not, let $I$ be a $B$-irreducible ideal of $\mathfrak{L}.$ By Lemma \ref{213.6}, $\mathfrak{L} = I \dot{+} I^{\bot}.$ Applying the induction hypothesis to $I$ and $I^{\bot}.$ This ends the proof of the proposition.
\end{proof}
\end{proposition}

Now we are ready to study central descending series $(C^{n}(\mathfrak{L}))_{n\geq 0}$ and central ascending series $(C_{n}(\mathfrak{L}))_{n\geq 0}.$ We will prove that $(C^{r}(\mathfrak{L}))^{\bot} = C_{r}(\mathfrak{L})$ for a quadratic Leibniz triple system $(\mathfrak{L},B).$ Also we conclude that the centralizer $Z_{\mathfrak{L}}(I)$ of a quadratic Leibniz triple system $(\mathfrak{L},B)$ contains $I^\bot$ if and only if $I$ is an ideal of $\mathfrak{L}.$ This work is inspired by {\rm\cite{WZ}}.

\begin{definition}
Let $\mathfrak{L}$ be a Leibniz triple system, the central descending series $(C^{n}(\mathfrak{L}))_{n\geq 0}$ is defined by $C^{0}(\mathfrak{L}):= \mathfrak{L}$ and $C^{n+1}(\mathfrak{L}):= \{C^{n}(\mathfrak{L}), \mathfrak{L}, \mathfrak{L}\} + \{\mathfrak{L}, C^{n}(\mathfrak{L}), \mathfrak{L}\} + \{\mathfrak{L}, \mathfrak{L}, C^{n}(\mathfrak{L})\}.$ And the central ascending series $(C_{n}(\mathfrak{L}))_{n\geq 0}$ is defined by $C_{0}(\mathfrak{L}):= 0$ and $C_{n+1}(\mathfrak{L}):= C(C_{n}(\mathfrak{L}))$ where $C(C_{r}(\mathfrak{L})):= \{x\in \mathfrak{L} ~|~\{x, \mathfrak{L}, \mathfrak{L}\} + \{\mathfrak{L}, x, \mathfrak{L}\} + \{\mathfrak{L}, \mathfrak{L}, x\} \subseteq C_{r}(\mathfrak{L})\}.$
\end{definition}
\begin{theorem}
Let $(\mathfrak{L}, B)$ be a quadratic Leibniz triple system. Then for all integer $r \geq 0,$ we have $(C^{r}(\mathfrak{L}))^{\bot} = C_{r}(\mathfrak{L}).$
\begin{proof}
We proceed by induction on $r.$

If $r = 0,$ $(C^{0}(\mathfrak{L}))^{\bot} = \mathfrak{L}^{\bot} = 0 = C_{0}(\mathfrak{L}),$ the equality holds.

Suppose the theorem is true for $r$ and we will show that it is true for $r+1.$ Let $w \in C_{r+1}(\mathfrak{L}).$ Note that $u \in C^{r+1}(\mathfrak{L})$ is a finite sum $u = \Sigma\Big(\{x_{i}, y_{i}, z_{i}\} + \{y_{i}^{\prime}, x_{i}^{\prime}, z_{i}^{\prime}\} + \{y_{i}^{\prime\prime}, z_{i}^{\prime\prime}, x_{i}^{\prime\prime}\}\Big),$ for all $x_{i}, x_{i}^{\prime}, x_{i}^{\prime\prime} \in C^{r}(\mathfrak{L})$ and $y_{i}, z_{i}, y_{i}^{\prime}, z_{i}^{\prime}, y_{i}^{\prime\prime}, z_{i}^{\prime\prime} \in \mathfrak{L}.$ Then
\begin{align*}
B(u, w) =& \Sigma B(\{x_{i}, y_{i}, z_{i}\} + \{y_{i}^{\prime}, x_{i}^{\prime}, z_{i}^{\prime}\} + \{y_{i}^{\prime\prime}, z_{i}^{\prime\prime}, x_{i}^{\prime\prime}\}, w)\\
        =& \Sigma B(\{x_{i}, y_{i}, z_{i}\}, w) + \Sigma B(\{y_{i}^{\prime}, x_{i}^{\prime}, z_{i}^{\prime}\}, w) + \Sigma B(\{y_{i}^{\prime\prime}, z_{i}^{\prime\prime}, x_{i}^{\prime\prime}\}, w)\\
        =& \Sigma B(x_{i}, \{w, z_{i}, y_{i}\}) + \Sigma B(x_{i}^{\prime}, \{z_{i}^{\prime}, w, y_{i}^{\prime}\}) + \Sigma B(x_{i}^{\prime\prime}, \{z_{i}^{\prime\prime}, y_{i}^{\prime\prime}, w\}).
\end{align*}
Since $(C^{r}(\mathfrak{L}))^{\bot} = C_{r}(\mathfrak{L})$ and
\begin{align*}
\{w, z_{i}, y_{i}\} &\in \{C_{r+1}(\mathfrak{L}), \mathfrak{L}, \mathfrak{L}\} \subseteq C_{r}(\mathfrak{L}),\\
\{z_{i}^{\prime}, w, y_{i}^{\prime}\} &\in \{\mathfrak{L}, C_{r+1}(\mathfrak{L}), \mathfrak{L}\} \subseteq C_{r}(\mathfrak{L}),\\
\{z_{i}^{\prime\prime}, y_{i}^{\prime\prime}, w\} &\in \{\mathfrak{L}, \mathfrak{L}, C_{r+1}(\mathfrak{L})\} \subseteq C_{r}(\mathfrak{L}),
\end{align*}
we have $B(x_{i}, \{w, z_{i}, y_{i}\}) = B(x_{i}^{\prime}, \{z_{i}^{\prime}, w, y_{i}^{\prime}\}) = B(x_{i}^{\prime\prime}, \{z_{i}^{\prime\prime}, y_{i}^{\prime\prime}, w\}) = 0,$ hence $B(u,w)=0.$ That means $C_{r+1}(\mathfrak{L}) \subseteq (C^{r+1}(\mathfrak{L}))^{\bot}.$

Conversely, for any $x \in (C^{r+1}(\mathfrak{L}))^{\bot}$ and for any $y, z \in \mathfrak{L},$ $w \in C^{r}(\mathfrak{L}),$ since $\{w, z, y\},$ $\{z, w, y\},$ $\{z, y, w\} \in C^{r+1}(\mathfrak{L}),$ we have
\begin{align*}
B(\{x, y, z\}, w) &= B(x, \{w, z, y\}) = 0,\\
B(\{y, x, z\}, w) &= B(x, \{z, w, y\}) = 0,\\
B(\{y, z, x\}, w) &= B(x, \{z, y, w\}) = 0,
\end{align*}
then $\{x, y, z\}, \{y, x, z\}, \{y, z, x\} \in (C^{r}(\mathfrak{L}))^{\bot} = C_{r}(\mathfrak{L}),$ i.e., $x \in C_{r+1}(\mathfrak{L}),$ that  means $(C^{r+1}(\mathfrak{L}))^{\bot} \subset C_{r+1}(\mathfrak{L}).$ This ends the proof of the theorem.
\end{proof}
\end{theorem}

\begin{definition}{\rm\cite{Ma}}\label{def254.9}
Let $\mathfrak{L}$ be a Leibniz triple system and $I$ a non-empty subset of $\mathfrak{L}.$ We call
\begin{align*}
Z_{\mathfrak{L}}(I)\! = \! \{x \in \mathfrak{L}~|~\{x, I, \mathfrak{L}\}=\{x, \mathfrak{L}, I\}=\{\mathfrak{L}, x, I\}=\{\mathfrak{L}, I, x\}=\{I, x, \mathfrak{L}\}=\{I, \mathfrak{L}, x\}=0\}
\end{align*}
the centralizer of $I$ in $\mathfrak{L}$.
\end{definition}

Combing Definition \ref{def254.9} and Lemma \ref{213.6}, one has
\begin{proposition}\label{213.11}
For a finite dimensional quadratic Leibniz triple system $(\mathfrak{L}, B),$ a subspace $I$ of $\mathfrak{L}$ is an ideal if and only if $Z_{\mathfrak{L}}(I)$ contains $I^{\bot}.$
\end{proposition}

\subsection{$T^{*}$-extension of a quadratic Leibniz triple system}

In this section, we give the method of $T^{*}$-extension for a Leibniz triple system $\mathfrak{L}.$ We will prove that every $2n$-dimensional Leibniz triple system $\mathfrak{L}$ is isomorphic to a $T^*$-extension of another Leibniz triple system if $\mathfrak{L}$ contains an $n$-dimensional isotropic ideal.

\begin{lemma}
Let $(\mathfrak{L}, \{\cdot, \cdot, \cdot\})$ be a Leibniz triple system and $\mathfrak{L}^{*}$ the dual space of $\mathfrak{L}.$ Consider a trilinear map $\theta : \mathfrak{L} \times \mathfrak{L} \times \mathfrak{L} \rightarrow \mathfrak{L}^{*}$ and the linear space $\mathfrak{L}\dot{+} \mathfrak{L}^{*}$ endowed with the triple product
\begin{align}
\{x+f, y+g, z+h\}_{\theta} = \{x, y, z\} + \theta(x, y, z) + \{f, y, z\}+ \{x, g, z\} + \{x, y, h\},\label{224.11}
\end{align}
for all $x, y, z \in \mathfrak{L}$ and $f, g, h \in \mathfrak{L}^{*},$ $(f,g,h)$ is a dual representation of $\mathfrak{L}.$ The linear space $\mathfrak{L}\dot{+} \mathfrak{L}^{*}$ endowed with $\{\cdot, \cdot, \cdot\}_{\theta}$ is a Leibniz triple system if and only if the map $\theta$ is a $3$-cocycle.
\begin{proof}
Example \ref{ex222.6} says that $\mathfrak{L}\dot{+} \mathfrak{L}^{*}$ is a Leibniz triple system for $\theta=0$. Then one observes that $\mathfrak{L}\dot{+} \mathfrak{L}^{*}$  in the assertion becomes a Leibniz triple system if and only if $\theta$ satisfies
\begin{gather}
\begin{aligned}
 &\theta(x,y,\{z,u,v\}) +  \{x,y,\theta(z,u,v)\}\\
  =& \theta(\{x,y,z\},u,v) - \theta(\{x,y,u\},z,v) - \theta(\{x,y,v\},z,u)+ \theta(\{x,y,v\},u,z)\\
 &+ \{\theta(x,y,z),u,v\} - \{\theta(x,y,u),z,v\}-\{\theta(x,y,v),z,u\}+ \{\theta(x,y,v),u,z\}, \label{224.8}
\end{aligned}
\end{gather}
and
\begin{gather}
\begin{aligned}
 &\theta(x,\{y,z,u\},v\})+\{x,\theta(y,z,u),v\}\\
   =& \theta(\{x,y,z\},u,v) - \theta(\{x,z,y\},u,v) - \theta(\{x,u,y\},z,v)+ \theta(\{x,u,z\},y,v)\\
 &+ \{\theta(x,y,z),u,v\} - \{\theta(x,z,y),u,v\}-\{\theta(x,u,y)z,v\}+ \{\theta(x,u,z),y,v\},\label{224.9}
\end{aligned}
\end{gather}
i.e., $\theta$ is a $3$-cocycle.
\end{proof}
\end{lemma}

\begin{definition}
Let $\mathfrak{L},$ $\mathfrak{L}^*,$ $\theta$ and $\{\cdot,\cdot,\cdot\}_{\theta}$ be as above. Then the Leibniz triple system $\mathfrak{L}\dot{+} \mathfrak{L}^{*}$ endowed with $\{\cdot,\cdot,\cdot\}_{\theta}$ is called a $T^*$-extension of $\mathfrak{L},$ denoted by  $T^*_\theta\mathfrak{L}.$ In particular, when $\theta=0,$ we denote it by $T^*\mathfrak{L}$ called the trivial $T^*$-extension of $\mathfrak{L}.$
\end{definition}

\begin{remark}\label{re4.13}
For any $T^*$-extension of $\mathfrak{L},$ it is easy to check that $\mathfrak{L}^*$ is an ideal of $T^*_\theta\mathfrak{L}$ and $\mathfrak{L}\simeq T^*_\theta\mathfrak{L}/\mathfrak{L}^*.$
\end{remark}

\begin{proposition}
Consider a bilinear form $B$ on $\mathfrak{L}\dot{+}\mathfrak{L}^{*}$ by
\begin{align}\label{224.55}
B(x+f, y+g) = g(x) + f(y),
\end{align}
for all $x, y\in \mathfrak{L},$ $f,g\in \mathfrak{L}^{*},$ if
\begin{align}\label{224.66}
\theta(x,y,z)(u)= \theta(y,x,u)(z) = \theta(u,z,y)(x),
\end{align}
for all $x,y,z,u \in \mathfrak{L},$ then $B$ is a non-degenerate symmetric invariant bilinear form on $\mathfrak{L}\dot{+} \mathfrak{L}^{*}.$ The pair $(\mathfrak{L}\dot{+} \mathfrak{L}^{*}, \{\cdot, \cdot, \cdot\}_{\theta}, B)$ is a quadratic Leibniz triple system.
\begin{proof}
It is sufficiency to prove the symmetric, non-degenerate and invariance of $B.$

For the left-invariance of $B,$ note that, for any $x, y, z, u \in \mathfrak{L}$ and $f, g, h, l \in \mathfrak{L}^{*},$
\begin{align*}
&B(L_{(x+f,y+g)}(z+h), u+l)= B(\{x+f,y+g, z+h\}_{\theta}, u+l)\\
                          =& B(\{x,y,z\}+ \theta(x,y,z)+\{f,y, z\} + \{x,g, z\} +\{x, y,h\}, u+l)\\
                          =& l(\{x,y,z\})+\theta(x,y,z)(u) + \{f,y, z\}(u) + \{x,g, z\}(u) + \{x, y,h\}(u)\\
                          =& l(\{x,y,z\}) + f(\{u,z,y\}) + g(\{z,u,x\}) + h(\{y,x,u\})+ \theta(x,y,z)(u),
\end{align*}
\begin{align*}
&B(z+h, L_{(y+g,x+f)}(u+l))= B(z+h, \{y+g, x+f, u+l\}_{\theta})\\
                          =& B(z+h, \{y,x,u\} + \theta(y,x,u)+\{g,x, u\} + \{y,f, u\} + \{y, x,l\})\\
                          =& h(\{y,x,u\}) +\{g,x, u\}(z) + \{y,f, u\}(z) + \{y, x,l\}(z)+\theta(y,x,u)(z)\\
                          =& h(\{y,x,u\}) + g(\{z,u,x\}) + f(\{u,z,y\}) + l(\{x,y,z\})+\theta(y,x,u)(z),
\end{align*}
which coincide Eq. (\ref{224.66}). Similarly, one can prove that right-invariance of $B.$ So we deduce the invariance of $B$.

As for the non-degenerate of $B$, let $x+f \in (\mathfrak{L}\dot{+} \mathfrak{L}^{*})^{\bot}.$ For any $y \in \mathfrak{L}$ we have $B(x+f, y+0)= f(y) +0(x) = 0,$ which implies $f=0.$
For any $g \in \mathfrak{L}^{*},$ $B(x+f, 0+g)= f(0) +g(x) = 0,$ we have $x=0,$ hence $B$ is non-degenerate.
\end{proof}
\end{proposition}

Now, we are ready to coincide that a quadratic Leibniz triple system $(\mathfrak{L}, B)$ and $T^*$-extension of a Leibniz triple system.

\begin{lemma}\label{lm224.15}
Let $(\mathfrak{L}, B)$ be a quadratic Leibniz triple system of  $2n$-dimension and $I$ an isotropic $({\rm i.e.}, I \subseteq I^{\bot})$ $n$-dimensional subspace of $\mathfrak{L}.$ Then $I$ is an ideal of $\mathfrak{L}$ if and only if $I$ is abelian ideal of $\mathfrak{L}$, that is $\{I,\mathfrak{L},\mathfrak{L}\}=\{\mathfrak{L},I,\mathfrak{L}\}=\{\mathfrak{L},\mathfrak{L},I\}=0.$
\begin{proof}
Note that dim$I$ + dim$I^\bot$ = $n$ + dim$I^\bot$ = $2n,$ dim$I^\bot =$ dim$I=n$ and $I\subset I^\bot,$ we have $I= I^\bot.$

If $I$ is an ideal of $\mathfrak{L},$ then by right-invariance of $B$ $$B(\mathfrak{L},\{I,I,\mathfrak{L}\})=B(\{\mathfrak{L},\mathfrak{L},I\},I)\subset B(I,I)=B(I,I^\bot)=0,$$ which implies $\{I,I,\mathfrak{L}\}\subset \mathfrak{L}^\bot=0.$ Similarly, $\{I,\mathfrak{L},I\}=\{\mathfrak{L},I,I\}=0.$

Conversely, it follows from  $$B(I,\{I,\mathfrak{L},\mathfrak{L}\})=B(\{\mathfrak{L},I,I\},I)=B(0,I)=0$$ that  $\{I,\mathfrak{L},\mathfrak{L}\}\subset I^\bot=I.$ Similarly, $\{\mathfrak{L},I,\mathfrak{L}\}\subset I$ and $\{\mathfrak{L},\mathfrak{L},I\}\subset I.$ This implies that $I$ is an ideal of $\mathfrak{L}.$
\end{proof}
\end{lemma}

\begin{theorem}
Let $(\mathfrak{L}, B)$ be a quadratic Leibniz triple system of $2n$-dimension. Then $(\mathfrak{L}, B)$ is isometric to a $T^{*}$-extension $(T^*_\theta J, B_J)$ if and only if $\mathfrak{L}$ contains an isotropic ideal $I$ of $n$-dimensional. Furthermore, $J$ is isomorphic to $\mathfrak{L}/I.$
\begin{proof}
The necessity is clear by Remark \ref{re4.13}, since $J^*$ is an ideal of $T^*_\theta J$ and $J\simeq T^*_\theta J/J^*$ and $B(J^*, J^*)=0,$ $J^*\subset (J^*)^\bot,$ i.e., $J^*$ is isotropic.

For the sufficiency, let $I$ be an $n$-dimensional isotropic ideal of $\mathfrak{L}.$ Suppose $J=\mathfrak{L}/I$ and $p :\mathfrak{L} \rightarrow J$ be the canonical projection. We will show that $\mathfrak{L}$ is some $T^*$-extension of $J$ by two steps.

\emph{Step 1, construct a $T^*$-extension $T^*_\theta J$ of $J.$}

First, by the non-degenerate of $B$ we can find an isotropic complement subspace $J_{0}$ to $I$ in $\mathfrak{L},$ i.e., $\mathfrak{L} = J_{0} \dot{+} I$ and $J_0 \subset J_0^\bot.$ Then dim$\mathfrak{L}-$dim$I=n$ and by Lemma \ref{lm224.15}. we have $J_0^\bot=J_0.$

Denote $p_0$ $($resp. $p_1)$ the projection $\mathfrak{L} \rightarrow J_0$ $($resp. $\mathfrak{L} \rightarrow I)$ and define a linear map
\begin{align*}
\mu^B:I &\rightarrow J^*\\
       i&\mapsto\mu^B(i):\mu^B(i)(p(x)):=B(i,x).
\end{align*}
It is easily seen that $\mu^B$ is well-defined, since for $p(x)=p(y)$ then $x-y\in I,$ hence $B(i,x-y)\in B(i,I)=0,$ so $B(i,x)=B(i,y),$ which implies $\mu^B$ is well-defined  and it is easily seen that $\mu^B$ is linear. If $\mu^B(i)=\mu^B(j)$ then $\mu^B(i)(p(x))-\mu^B(j)(p(x))=B(i,x)-B(j,x),$ for all $x\in \mathfrak{L},$ which implies $i-j\in \mathfrak{L}^\bot=0,$ hence $\mu^B$ is injective. As dim$I=$ dim$J^*,$ $\mu^B$ is also surjective.

In addition, $\mu^B$ has the following properties:
\begin{align*}
\mu^B(\{a_1,a_2,i\})(pa_3)=&B(\{a_1,a_2,i\},a_3) = B(i,\{a_2,a_1,a_3\})\\ =&\mu^B(i)(\{pa_2,pa_1,pa_3\})=\{pa_1,pa_2,\mu^B(i)\}(pa_3),
\end{align*}
i.e., $\mu^B(\{a_1,a_2,i\})=\{pa_1,pa_2,\mu^B(i)\},$
for any $a_1,a_2,a_3\in \mathfrak{L}$ and $i\in I.$ Further, a similar computation shows that,
\begin{align*}
\mu^B(\{i,a_1,a_2\})=&\{\mu^B(i),pa_1,pa_2\}\\
\mu^B(\{a_1,i,a_2\})=&\{pa_1,\mu^B(i),pa_2\}.
\end{align*}

Define a trilinear map
\begin{align*}
\theta : J\times J \times J &\rightarrow J^*\\
      (p(b_1),p(b_2),p(b_3))&\mapsto \mu^B(p_1(\{b_1,b_2,b_3\})),
\end{align*}
for all $b_1,b_2,b_3\in J_0.$ Then $\theta$ is well-defined, since $p|_{J_0}:J_0\rightarrow J_0/I\simeq \mathfrak{L}/I=J$ is a linear isomorphism and $\theta\in C^3(J,J^*)$ it is easy to check that $\theta\in Z^3(J,J^*).$ Now, define the bracket on $J\dot{+} J^*$ by (\ref{224.11}). Then $J\dot{+} J^*$ is a Leibniz triple system.

\emph{Step 2, construct a isomorphism from $\mathfrak{L}$ to $J\dot{+} J^*.$}

Let $\phi: \mathfrak{L} \rightarrow J\dot{+} J^*$ be a  linear map, via $\phi(b_0 + i)= pb_0 + \mu^B(i),$ for any $b_0 \in J_0, i\in I.$ Since the restriction of $p$ to $J_0$ and $\mu^B$ are linear isomorphisms, the map $\phi$ is also a linear isomorphism. Note that for any $b_1,b_2, b_3 \in J_0$ and $i_1,i_2,i_3 \in I,$ we have
\begin{align*}
&\phi(\{b_1 + i_1,b_2 + i_2,b_3 + i_3\})\\
=&\phi\Big (\{b_1,  b_2, b_3\}+\{b_1,  b_2, i_3\} + \{i_1,  b_2, b_3\} + \{b_1,  i_2, b_3\}\Big)\\
=& \phi\Big (p_0(\{b_1,  b_2, b_3\}) +p_1(\{b_1,  b_2, b_3\}) + \{b_1,  b_2, i_3\} + \{i_1,  b_2, b_3\} + \{b_1,  i_2, b_3\}\Big)\\
=&p(p_0(\{b_1,  b_2, b_3\})) + \mu^B\Big(p_1(\{b_1,  b_2, b_3\})+ \{b_1,  b_2, i_3\} + \{i_1,  b_2, b_3\} + \{b_1,  i_2, b_3\}\Big)\\
=&p(\{b_1,  b_2, b_3\}) + \Omega (pb_1,  pb_2, pb_3) + \{pb_1,  pb_2, \mu^B(i_3)\} + \{\mu^B(i_1),  pb_2, pb_3\} + \{pb_1,  \mu^B(i_2), pb_3\}\\
=&\{pb_1+\mu^B(i_1),pb_2+\mu^B(i_2),pb_3+\mu^B(i_3)\}\\
=&\{\phi(b_1 + i_1), \phi(b_2 + i_2), \phi(b_3 + i_3)\},
\end{align*}
where we have made use of the definition of $\phi$ and $\theta,$ and the above properties of $\mu^B.$ Then $\phi$ is a Leibniz triple system isomorphism. Furthermore, we have
\begin{align*}
B_J(\phi(x+i),\phi(x^\prime+i^\prime))=&B_J(p(x)+\mu^B(i),p(x^\prime)+\mu^B(i^\prime))\\
                                      =&\mu^B(i^\prime)(p(x))+\mu^B(i)(p(x^\prime))\\
                                      =&B(i^\prime,x)+B(i,x^\prime)\\
                                      =&B(x+i,x^\prime+i^\prime).
\end{align*}
Then $\phi$ is isometric. We also have $B_J$ is non-degenerate following from the non-degenerate property of $B.$ So $(J\dot{+} J^*, B_J)$ is a quadratic Leibniz triple system. In this way, we get a $T^*$-extension $T^*_\theta J$ of $J,$ and consequently, $(\mathfrak{L},B)$ and $(T^*_\theta J, B_J)$ are isometric as required.
\end{proof}
\end{theorem}

\subsection{Symplectic form of a Leibniz triple system}
In this subsection, we will introduce the symplectic form of a Leibniz triple system, and we give a necessary and sufficient condition for a quadratic Leibniz triple system to admit a symplectic form.

\begin{definition}
Let $(\mathfrak{L}, \{\cdot, \cdot, \cdot\})$ be a Leibniz triple system. A symplectic form  is a non-degenerate skewsymmetric bilinear form $\omega$ satisfying
\begin{align}
\omega(u, \{x,y,z\}) + \omega(x, \{u,z,y\})+ \omega(y, \{z,u,x\}) + \omega(z, \{y,x,u\}) = 0,\label{****}
\end{align}
for all $x, y, z, u\in \mathfrak{L}.$ A Leibniz triple system endowed with a symplectic form is called symplectic Leibniz triple system, denoted by $(\mathfrak{L}, \omega).$
\end{definition}

\begin{definition}
A symmetric $($resp. skewsymmetric$)$ linear map $f$ on $\mathfrak{L}$ with respect to $B$ satisfying $B(f(x), y) = B(x, f(y))$ $({\rm resp}. ~B(f(x), y) = -B(x, f(y))),$ for all $x, y \in \mathfrak{L}.$
\end{definition}

\begin{theorem}\label{214.5}
Let $(\mathfrak{L},B)$ be a quadratic Leibniz triple system. There is a symplectic form $\omega$ on $\mathfrak{L}$ if and only if there is an invertible derivation of $\mathfrak{L}$ which is skewsymmetric with respect to $B.$
\begin{proof}
For the sufficiency, let $D$ be an invertible and skewsymmetric derivation on $\mathfrak{L}$ with respect to $B.$ Define $\omega(x,y) = B(D(x), y).$ In the following we prove that $\omega$ is skewsymmetric, non-degenerate and satisfies Eq. (\ref{****}). Since
$$\omega(x,y) = B(D(x), y) = -B(x, D(y)) =-B(D(y), x)=-\omega(y,x),$$
and $D$ is invertible, $B$ is non-degenerate, we have $\omega$ is skewsymmetric and non-degenerate. Furthermore, by the skewsymmetric of $D$ with respect to $B$ and the invariance of $B,$
\begin{align*}
\omega(u, \{x, y, z\}) =& B(D(u), \{x, y, z\})\\
                       =& -B(u, D(\{x, y, z\}))\\
                       =& -B(u, \{D(x), y, z\}) - B(u, \{x, D(y), z\}) - B(u, \{x, y, D(z)\})\\
                       =& -B(\{u, z, y\}, D(x)) - B(\{z, u, x\}, D(y)) - B(\{y, x,u\}, D(z))\\
                       =& -B(D(x), \{u, z, y\}) - B(D(y), \{z, u, x\}) - B(D(z), \{y, x,u\})\\
                       =& -\omega(x, \{u, z, y\}) - \omega(y, \{z, u, x\}) - \omega(z, \{y, x,u\}),
\end{align*}
then we have
\begin{align*}
\omega(u, \{x, y, z\})+\omega(x, \{u, z, y\}) + \omega(y, \{z, u, x\}) + \omega(z, \{y, x,u\}) =0,
\end{align*}
for all $x, y, z, u \in \mathfrak{L}.$ So $\omega$ is a symplectic form.

For the necessity, if $\omega$ is a symplectic form on $\mathfrak{L}.$ Define $D$ by $\omega(x,y) = B(D(x), y),$ for all $x,y\in \mathfrak{L},$ which is well-defined since the non-degenerate of $B$ and $B,$ $\omega$ are bilinear form deduce that $D$ is a linear map. By the definition of $\omega$ we have
\begin{align*}
B(D(u), \{x, y, z\})+B(D(x), \{u, z, y\}) + B(D(y), \{z, u, x\}) + B(D(z), \{y, x, u\})=0.
\end{align*}
Since $D$ is skewsymmetric, $B$ is invariant and symmetric,
\begin{align*}
-B(D(\{x, y, z\}),u)+ B(\{D(x), y, z\}, u\}) + B(\{x, D(y), z\}, u\}) + B(\{x, y, D(z)\}, u) = 0,
\end{align*}
i.e.,
\begin{align*}
B(D(\{x, y, z\}), u)=B(\{D(x), y, z\}, u\}) + B(\{x, D(y), z\}, u\}) + B(\{x, y, D(z)\}, u).
\end{align*}

Again $B$ is non-degenerate $D$ becomes a derivation on $\mathfrak{L}.$ $D$ is skewsymmetric, since
\begin{align*}
B(D(x),y)= \omega(x,y)=-\omega(y,x)= -B(D(y),x)=-B(x,D(y)).
\end{align*}

For the invertible of $D,$ define $D^{\prime}$ by $\omega(D^{\prime}(x),y)=B(x,y)$ for any $x,y\in \mathfrak{L},$ which is well-defined since the non-degenerate of $B,$  note that
\begin{align*}
&\omega(D^{\prime}D(x),y)= B(D(x),y)=\omega(x,y),\\
&B(DD^{\prime}(x),y) = \omega(D^{\prime}(x),y) = B(x,y),
\end{align*}
it follows that $D^{\prime}D= DD^{\prime}=$id. Thus $D$ is invertible. So we have finished the proof of this theorem.
\end{proof}
\end{theorem}

\begin{remark}
The skewsymmetric derivation $D$ of a quadratic Leibniz triple system $(\mathfrak{L},B)$ in the above theorem is also skewsymmetric with respect to the symplectic form $\omega.$ In fact, for any $x,y\in \mathfrak{L},$
\begin{align*}
\omega(D(x),y)= B(D^{2}(x),y) = -B(D(x),D(y))=-\omega(x,D(y)).
\end{align*}
\end{remark}

\begin{corollary}\label{coro4.21}
Let $T^*\mathfrak{L}$ be a $T^{*}$-extension of a Leibniz triple system $\mathfrak{L}$ admitting a non-degenerate symmetric invariant bilinear form $B$ defined by Eq. $(\ref{224.55})$, if $T^*\mathfrak{L}$ has an invertible skewsymmetric derivation $D,$ then there exists a symplectic form for $T^*\mathfrak{L}.$
\end{corollary}

By the conclusion of Lemma \ref{coro4.21}, we give the following example for a symplectic Leibniz triple system.
\begin{example}
Let $(\mathfrak{L}, \{\cdot, \cdot, \cdot\}_{\mathfrak{L}})$ be a Leibniz triple system. Consider the $G_{n} = t\mathbb{F}[t]/t^{n}\mathbb{F}[t],$ where $n$ is a nonzero positive integer. Define the product on $\mathfrak{L}_{n} = \mathfrak{L} \otimes G_{n},$ as follows,
\begin{align*}
\{x\otimes \overline{t^{p}}, y\otimes \overline{t^{q}}, z\otimes \overline{t^{r}}\} = \{x, y, z\}_{\mathfrak{L}} \otimes \overline{t^{p+q+r}},
\end{align*}
where $x, y, z\in \mathfrak{L},$ $p, q, r\in \mathbb{N}/\{0\}.$ Then $\mathfrak{L}_{n}$ is a nilpotent Leibniz triple system.  By the universal property of tensor product define a linear map $D$ on $\mathfrak{L}_{n}$ by $D(x\otimes \overline{t^{p}}) = p(x\otimes \overline{t^{p}}),$ for any $x \in \mathfrak{L},$ and $p\in \{1, ..., n-1\}$ it is easy to check that $D$ is an invertible derivation of $\mathfrak{L}_{n}.$

Consider the trivial $T^*$-extension deduce the quadratic Leibniz triple system  on $\mathfrak{L}_{n},$ which is denoted by $\widetilde{\mathfrak{L}_{n}} := \mathfrak{L}_{n} \dot{+} \mathfrak{L}_{n}^{*},$ when $\widetilde{\mathfrak{L}_{n}}$ endowed with the derivation $\widetilde{D}(x+f) = D(x) - f\circ D,$ for $x\in \mathfrak{L}_{n},$ $f\in \mathfrak{L}_{n}^{*},$ then $\widetilde{D}$ is invertible with $\widetilde{D}^{-1}=D^{-1}(x)-f\circ D^{-1},$ so $\widetilde{D}$ is an invertible derivation on $\widetilde{\mathfrak{L}_{n}}$ which is skewsymmetric with respect to $B,$ since $B(\widetilde{D}(x+f),y+h)=h\circ D(x)-f\circ D(y),$ hence the quadratic Leibniz triple system $(\widetilde{\mathfrak{L}_{n}},B)$ admits a symplectic structure.
\end{example}

\section{$1$-parameter formal deformation of Leibniz triple systems}
At last, we use the cohomology theory defined in Section 2 to study the $1$-parameter formal deformation of Leibniz triple systems. Let $(\mathfrak{L},\{\cdot,\cdot,\cdot\})$ be a Leibniz triple system and $\mathbb{F}[[t]]$ be the formal power series over $\mathbb{F}.$ Suppose that $\mathfrak{L}[[t]]$ is the set of formal power series over $\mathfrak{L}.$ Then for an $\mathbb{F}$-trilinear map $f:\mathfrak{L}\times \mathfrak{L}\times \mathfrak{L}\rightarrow \mathfrak{L},$ it is natural to extend it to an $\mathbb{F}[[t]]$-trilinear map $f:\mathfrak{L}[[t]]\times \mathfrak{L}[[t]]\times \mathfrak{L}[[t]]\rightarrow \mathfrak{L}[[t]],$ by
\begin{align*}
f\Big(\sum\limits_{i\geq0}^{}x_it^i,\sum\limits_{j\geq0}^{}y_jt^j,\sum\limits_{k\geq0}^{}z_kt^k\Big)
=\sum\limits_{i,j,k\geq0}^{}f(x_i,y_j,z_k)t^{i+j+k}.
\end{align*}

\begin{definition}
Let $(\mathfrak{L},\{\cdot,\cdot,\cdot\})$ be a Leibniz triple system over $\mathbb{F}.$ A $1$-parameter formal deformation of $(\mathfrak{L},\{\cdot,\cdot,\cdot\})$ is a formal power series $d_t : \mathfrak{L}[[t]]\times \mathfrak{L}[[t]]\times \mathfrak{L}[[t]]\rightarrow \mathfrak{L}[[t]]$ of the form $d_t(x,y,z)=\sum\limits_{i\geq 0}^{}d_i(x,y,z)t^i$ where each $d_i$ is an $\mathbb{F}$-trilinear map $d_i:\mathfrak{L}\times \mathfrak{L}\times \mathfrak{L}\rightarrow \mathfrak{L}$ $($extended to be $\mathbb{F}[[t]]$-trilinear$)$ and $d_0(x,y,z)=\{x,y,z\},$ such that the following identities hold
\begin{gather}
\begin{aligned}
d_t(u,v,d_t(x,y,z))=&d_t(d_t(u,v,x),y,z)-d_t(d_t(u,v,y),x,z)-d_t(d_t(u,v,z),x,y)\\
                    &+d_t(d_t(u,v,z),y,x)\label{226.1},\\
\end{aligned}
\end{gather}
\begin{gather}
\begin{aligned}
d_t(u,d_t(x,y,z),v)=&d_t(d_t(u,x,y),z,v)-d_t(d_t(u,y,x),z,v)-d_t(d_t(u,z,x),y,v)\\
                    &+d_t(d_t(u,z,y),x,v)\label{226.2}.
\end{aligned}
\end{gather}
Eqs.$(\ref{226.1})$ and $(\ref{226.2})$ are called the deformation equations and give a Leibniz triple system structure on $\mathfrak{L}_t=(\mathfrak{L}[[t]],d_t).$
\end{definition}

Note that Eqs. $(\ref{226.1})$ and $(\ref{226.2})$ can be expressed as
\begin{gather}
\begin{aligned}
\sum\limits_{i,j\geq0}^{}d_i(u,v,d_j(x,y,z))=&\sum\limits_{i,j\geq0}^{}d_i(d_j(u,v,x),y,z)
-\sum\limits_{i,j\geq0}^{}d_i(d_j(u,v,y),x,z)\\
&-\sum\limits_{i,j\geq0}^{}d_i(d_j(u,v,z),x,y)+\sum\limits_{i,j\geq0}^{}d_i(d_j(u,v,z),y,x),\label{226.3}
\end{aligned}
\end{gather}
\begin{gather}
\begin{aligned}
\sum\limits_{i,j\geq0}^{}d_i(u,d_j(x,y,z),v)=&\sum\limits_{i,j\geq0}^{}d_i(d_j(u,x,y),z,v)
-\sum\limits_{i,j\geq0}^{}d_i(d_j(u,y,x),z,v)\\
&-\sum\limits_{i,j\geq0}^{}d_i(d_j(u,z,x),y,v)+\sum\limits_{i,j\geq0}^{}d_i(d_j(u,z,y),x,v).\label{226.4}
\end{aligned}
\end{gather}
Then
\begin{gather}
\begin{aligned}
\sum\limits_{i+j=r}^{}&\Big(d_i(d_j(u,v,x),y,z)-d_i(d_j(u,v,y),x,z)-d_i(d_j(u,v,z),x,y)+d_i(d_j(u,v,z),y,x)\\
&-d_i(u,v,d_j(x,y,z))\Big)=0,\label{226.5}
\end{aligned}
\end{gather}
\begin{gather}
\begin{aligned}
\sum\limits_{i+j=r}^{}&\Big(d_i(d_j(u,x,y),z,v)-d_i(d_j(u,y,x),z,v)-d_i(d_j(u,z,x),y,v)+d_i(d_j(u,z,y),x,v)\\
&-d_i(u,d_j(x,y,z),v)\Big)=0.\label{226.6}
\end{aligned}
\end{gather}
\begin{definition}
Suppose that $d_t(x,y,z)=\sum\limits_{i\geq 0}^{}d_i(x,y,z)t^i$ is a $1$-parameter formal deformation of $(\mathfrak{L},\{\cdot,\cdot,\cdot\}).$ Then $d_1$ is called an infinitesimal deformation. In general, if $d_i=0,$ for $1 \leq i\leq n-1,$ $d_n\neq 0,$ then $d_n$ is called an $n$-infinitesimal deformation.
\end{definition}

\begin{proposition}\label{prop255.3}
The infinitesimal deformation $d_1$ is a $3$-cocycle in $C^3(\mathfrak{L},\mathfrak{L}).$ In general, $n$ infinitesimal deformation $d_n$ is also a $3$-cocycle in $C^3(\mathfrak{L},\mathfrak{L}).$
\begin{proof}
Taking $r=1$ in Eqs. (\ref{226.5}) and (\ref{226.6}), we obtain two equations for $d_1:$
\begin{gather}
\begin{aligned}
&\{d_1(u,v,x),y,z\}-\{d_1(u,v,y),x,z\}-\{d_1(u,v,z),x,y\}+\{d_1(u,v,z),y,x\}\\
&-\{u,v,d_1(x,y,z)\}+d_1(\{u,v,x\},y,z)-d_1(\{u,v,y\},x,z)-d_1(\{u,v,z\},x,y)\\
&+d_1(\{u,v,z\},y,x)-d_1(u,v,\{x,y,z\})=0,\label{226.7}
\end{aligned}
\end{gather}
\begin{gather}
\begin{aligned}
&\{d_1(u,x,y),z,v\})-\{d_1(u,y,x),z,v\}-\{d_1(u,z,x),y,v\}+\{d_1(u,z,y),x,v\}\\
&-\{u,d_1(x,y,z),v\}+d_1(\{u,x,y\},z,v)-d_1(\{u,y,x\},z,v)-d_1(\{u,z,x\},y,v)\\
&+d_1(\{u,z,y\},x,v)-d_1(u,\{x,y,z\},v)=0.\label{226.8}
\end{aligned}
\end{gather}
So $d_1$ is a $3$-cocycle with respect to adjoint module $\mathfrak{L}.$ Similarly, one can check that $d_n$ is a $3$-cocycle.
\end{proof}
\end{proposition}

Now, we characterize the $1$-parameter formal deformation through cohomology of $\mathfrak{L}.$

\begin{definition}
Let $(\mathfrak{L},\{\cdot,\cdot,\cdot\})$ be a Leibniz triple system. Suppose that $d_t(x,y,z)=\sum\limits_{i\geq0}^{} d_i(x,y,z)t^i$ and $d_t^\prime(x,y,z)=\sum\limits_{i\geq0}^{} d_i^\prime(x,y,z)t^i$ are two $1$-parameter formal deformations of $(\mathfrak{L},\{\cdot,\cdot,\cdot\}).$ They are called equivalent, denote by $d_t\sim d_t^\prime,$ if there is a formal isomorphism of $\mathbb{F}[[t]]$-modules
\begin{align*}
\phi_t(x)=\sum\limits_{i\geq0}^{}\phi_i(x)t^i:(\mathfrak{L}[[t]],d_t)\rightarrow (\mathfrak{L}[[t]],d_t^\prime),
\end{align*}
where each $\phi_i : \mathfrak{L}\rightarrow \mathfrak{L}$ is an $\mathbb{F}$-linear map $($extended to be $\mathbb{F}[[t]]$-linear$)$ and $\phi_0=$id$_\mathfrak{L},$ satisfying
\begin{align*}
\phi_td_t(x,y,z)=d_t^\prime (\phi_t(x),\phi_t(y),\phi_t(z)).
\end{align*}
If $d_1=d_2=... =0,$ $d_t=d_0$ is said to be the null deformation. A $1$-parameter formal deformation $d_t$ is called trivial if $d_t\sim d_0.$ A Leibniz triple system $(\mathfrak{L},\{\cdot,\cdot,\cdot\})$ is called analytically rigid, if every $1$-parameter formal deformation $d_t$ is trivial.
\end{definition}

\begin{theorem}
Let $d_t(x,y,z)=\sum\limits_{i\geq0}^{} d_i(x,y,z)t^i$ and $d_t^\prime(x,y,z)=\sum\limits_{i\geq0}^{} d_i^\prime(x,y,z)t^i$ be equivalent $1$-parameter formal deformations of $(\mathfrak{L},\{\cdot,\cdot,\cdot\}).$ Then $d_1$ and $d_1^\prime$ belong to the same cohomology class in $H^3(\mathfrak{L},\mathfrak{L}).$
\begin{proof}
Suppose that $\phi_t(x)=\sum\limits_{i\geq0}^{}\phi_i(x)t^i$ is the formal $\mathbb{F}[[t]]$-module isomorphism such that
\begin{align*}
\sum\limits_{i\geq0}^{}\phi_i\Big(\sum\limits_{j\geq0}^{}d_j(x,y,z)\Big)t^{i+j}=\sum\limits_{i\geq0}^{}d_i^\prime \Big(\sum\limits_{k\geq0}^{}\phi_k(x)t^k,\sum\limits_{l\geq0}^{}\phi_l(y)t^l,\sum\limits_{m\geq0}^{}\phi_m(z)t^m\Big)t^i.
\end{align*}
It follows that
\begin{align*}
\sum\limits_{i+j=n}^{}\phi_i(d_j(x,y,z))t^{i+j}=\sum\limits_{i+k+l+m=n}^{}d_i^\prime(\phi_k(x),\phi_l(y),\phi_m(z))t^{i+k+l+m}.
\end{align*}
In particular,
\begin{align*}
\sum\limits_{i+j=1}^{}\phi_i(d_j(x,y,z))=\sum\limits_{i+k+l+m=1}^{}d_i^\prime(\phi_k(x),\phi_l(y),\phi_m(z)),
\end{align*}
that is,
\begin{align*}
d_1(x,y,z)+\phi_1(\{x,y,z\})=\{\phi_1(x),y,z\} + \{x,\phi_1(y),z\} + \{x,y,\phi_1(z)\}+d_1^\prime(x,y,z),
\end{align*}
we have
\begin{align*}
d_1(x,y,z)-d_1^\prime(x,y,z)=\{\phi_1(x),y,z\} + \{x,\phi_1(y),z\} + \{x,y,\phi_1(z)\}-\phi_1(\{x,y,z\}).
\end{align*}
Then $d_1-d_1^\prime \in B^3(\mathfrak{L},\mathfrak{L}).$
\end{proof}
\end{theorem}

\begin{theorem}
Suppose that $(\mathfrak{L},\{\cdot,\cdot,\cdot\})$ is a Leibniz triple system such that $H^3(\mathfrak{L},\mathfrak{L})=0.$ Then $(\mathfrak{L},\{\cdot,\cdot,\cdot\})$ is analytically rigid.
\begin{proof}
Let $d_t$ be a $1$-parameter formal deformation of $(\mathfrak{L},\{\cdot,\cdot,\cdot\}).$ Suppose that $d_t=d_0+\sum\limits_{i\geq n}^{}d_it^i,$ and by Proposition \ref{prop255.3} $d_n\in Z^3(\mathfrak{L},V)=B^3(\mathfrak{L},\mathfrak{L}).$
It follows that there exists $f_n\in C^1(\mathfrak{L},\mathfrak{L}),$ such that
\begin{align}
d_n(x,y,z)=\{f_n(x),y,z\} + \{x,f_n(y),z\} + \{x,y,f_n(z)\}-f_n(\{x,y,z\}).\label{ddd}
\end{align}

Let $\phi_t=$id$_\mathfrak{L}-f_nt^n$ and $d_t^\prime(x,y,z)=\phi_t^{-1}d_t(\phi_t(x),\phi_t(y),\phi_t(z)).$ It is straightforward to prove that $d_t^\prime$ is a $1$-parameter formal deformation of $(\mathfrak{L},\{\cdot,\cdot,\cdot\})$ and $d_t\sim d_t^\prime.$ Suppose that $d_t^\prime=\sum\limits_{i\geq0}^{}d_i^\prime t^i,$ we express the coefficient in front of the $d_t^\prime$ in terms of $d_i^\prime.$  Then
\begin{align*}
({\rm id}_\mathfrak{L}-f_nt^n)\Big(\sum\limits_{i\geq 0}^{}d_i^\prime(x,y,z)t^i\Big)=\Big(d_0+\sum\limits_{i\geq n}^{}d_it^i\Big)(x-f_n(x)t^n,y-f_n(y)t^n,z-f_n(z)t^n),
\end{align*}
i.e.,
\begin{align*}
&\sum\limits_{i\geq 0}^{}d_i^\prime(x,y,z)t^i-\sum\limits_{i\geq 0}^{}f_nd_i^\prime(x,y,z)t^{i+n}\\
=&\{x,y,z\}-(\{f_n(x),y,z\}+\{x,f_n(y),z\}+\{x,y,f_n(z)\})t^n+(\{f_n(x),f_n(y),z\}\\
&+\{x,f_n(y),f_n(z)\}+\{f_n(x),y,f_n(z)\})t^{2n}-\{f_n(x),f_n(y),f_n(z)\}t^{3n}+\sum\limits_{i\geq n}^{}d_i(x,y,z)t^i\\
&-\sum\limits_{i\geq n}^{}(d_i(f_n(x),y,z)+d_i(x,f_n(y),z)+d_i(x,y,f_n(z)))t^{i+n}\\
&+\sum\limits_{i\geq n}^{}(d_i(f_n(x),f_n(y),z)+d_i(x,f_n(y),f_n(z))+d_i(f_n(x),y,f_n(z)))t^{i+2n}\\
&-\sum\limits_{i\geq n}^{}(d_i(f_n(x),f_n(y),f_n(z))t^{i+3n}.
\end{align*}
Then we have $d_1^\prime=... =d_{n-1}^\prime=0$ and
\begin{align*}
d_n^\prime(x,y,z)-f_n(\{x,y,z\})=-(\{f_n(x),y,z\} + \{x,f_n(y),z\} + \{x,y,f_n(z)\})+d_n(x,y,z).
\end{align*}
By Eq. (\ref{ddd}), $d_n^\prime=0$ and $d_t^\prime=d_0+\sum\limits_{i\geq n+1}^{}d_i^\prime t^i.$ Repeat this procedure, one observe that there is no nontrivial deformation of $(\mathfrak{L},\{\cdot,\cdot,\cdot\}),$ i.e., $(\mathfrak{L},\{\cdot,\cdot,\cdot\})$ is analytically rigid.
\end{proof}
\end{theorem}

At last, we define two $5$-cochains \rm Ob$_{n+1}^i(i=1,2)$ by
\begin{align*}
{\rm Ob}_{n+1}^1(u,v,x,y,z)=&\sum \limits_{i+j=n+1\atop i,j>0}^{}\Big(d_i(d_j(u,v,x),y,z)-d_i(d_j(u,v,y),x,z)\\
                     &-d_i(d_j(u,v,z),x,y)+d_i(d_j(u,v,z),y,x)-d_i(u,v,d_j(x,y,z))\Big),
\end{align*}
\begin{align*}
{\rm Ob}_{n+1}^2(u,v,x,y,z)=&\sum \limits_{i+j=n+1\atop i,j>0}^{}\Big(d_i(d_j(u,v,x),y,z)-d_i(d_j(u,x,v),y,z)\\
                     &-d_i(d_j(u,y,v),x,z)+d_i(d_j(u,y,x),v,z)-d_i(u,d_j(v,x,y),z)\Big).
\end{align*}

\begin{definition}
The $5$-cochains \rm Ob$_{n+1}^i(\mathfrak{L})\in C^5(\mathfrak{L},\mathfrak{L})(i=1,2)$ are called $(n+1)$-th obstruction cochain for extending a given deformation of order $n$ (i.e. $d_n=\sum\limits_{i=0}^{n}d^it^i$) to a deformation of $\mathfrak{L}$ of order $(n+1).$
\end{definition}

\begin{theorem}
Let $d_t$ be a deformation of $\mathfrak{L}$ of order $n.$ Then $d_t$ extends to a deformation of order $n+1,$ $d_{n+1}=\sum\limits_{i=0}^{n+1}d^it^i$ if and only if
\rm Ob$_{n+1}^i(\mathfrak{L})=\delta^3_id_{n+1}(i=1,2).$
\begin{proof}
For the necessity, it is obvious by the definition of \rm Ob$_{n+1}^i(\mathfrak{L})$ and Eqs. (\ref{226.5}) and (\ref{226.6}).

For the sufficiency, suppose \rm Ob$_{n+1}^i(\mathfrak{L})=\delta^3_id_{n+1}(i=1,2)$, for some $3$-cochain $d_{n+1}\in C^3(\mathfrak{L},\mathfrak{L}),$ that is \rm Ob$_{n+1}^i(\mathfrak{L})$ is a coboundary and  $\tilde{d_{t}}=d_t + d_{n+1}t^{n+1}.$ Observe that $\tilde{d_{t}}$ satisfies Eqs. (\ref{226.5}) and (\ref{226.6}) for $0\leq r\leq n+1$. So  $\tilde{d_{t}}$ is an extension of $d_t$ and is of order $n+1.$
\end{proof}
\end{theorem}

\end{document}